\documentclass[11pt]{article}

\usepackage[margin=1.0in]{geometry}
\usepackage{amsmath}
\usepackage{mathtools}
\usepackage{amsthm}
\usepackage{amsfonts}
\usepackage{fancyhdr}
\usepackage{graphicx}
\usepackage{float}
\usepackage{booktabs}
\usepackage{cmbright}
\usepackage{hyperref}

\newtheorem{theorem}{Theorem}[section]

\newcommand{\diff}{\textnormal{d}}

\newcommand{\bff}{\mathbf{f}}
\newcommand{\bfg}{\mathbf{g}}
\newcommand{\bfr}{\mathbf{r}}
\newcommand{\bfw}{\mathbf{w}}
\newcommand{\bfx}{\mathbf{x}}
\newcommand{\bfy}{\mathbf{y}}
\newcommand{\bfz}{\mathbf{z}}

\newcommand{\bfG}{\mathbf{G}}
\newcommand{\bfR}{\mathbf{R}}
\newcommand{\bfU}{\mathbf{U}}

\newcommand{\bbZ}{\mathbb{Z}}
\newcommand{\bbR}{\mathbb{R}}
\newcommand{\bbC}{\mathbb{C}}

\allowdisplaybreaks

\title{Learning Closure of Dynamical Systems with Kernel Ridge Regression}
\date{30 September 2026}

\author{
  Evan Habbershaw\thanks{
    Department of Mathematics, The Pennsylvania State University; Corresponding author: \texttt{ezh302@psu.edu}
  },
  John Harlim\thanks{
    Department of Mathematics and Institute for Computational and Data Sciences, The Pennsylvania State University
  }, and
  Senwei Liang\thanks{
    Department of Mathematics and Statistics, Texas Tech University
  }
}

\begin{document}

\pagenumbering{arabic}

\maketitle

\begin{abstract}

    \noindent
    We develop a closure modeling framework for identifying missing components of dynamical systems using Kernel Ridge Regression (KRR). The framework addresses two classes of closure problems: difference-equation closures arising in ODE and PDE settings, and algebraic closures arising from moment closure in kinetic equations. For the first class, we derive an error bound in an ODE setting that quantifies contributions from time integration, approximation of unresolved scales, and interpolation required to couple unresolved-scale effects to the resolved solver. Numerical experiments on the Lorenz-63 system and the Kuramoto–Sivashinsky equation demonstrate accurate long-horizon predictions and substantial improvements over an LSTM-based closure model. For the second class, we consider moment closure for a one-dimensional kinetic equation by modeling discrepancies between kinetic and macroscopic fluxes as a function of the resolved macroscopic variables. We compare global KRR models based on PCA coordinates with spatially local models. While the global model performs well for unimodal initial conditions, its accuracy deteriorates for bimodal initial conditions. Spatially local models with appropriate modeling inputs improve robustness and achieve higher predictive accuracy.

\end{abstract}


\section{Introduction}

The closure of unresolved components in dynamical systems is a longstanding problem that has been studied under various names, including model error, subgrid-scale parameterization, and closure modeling \cite{crommelin2008subgrid,bh:14,hmm:14,kwasniok2012data,lu2016comparison,lu2017accounting,mh:13,harlim_2021,levine2022framework}.
In this work, we consider two classes of closure problems.
The first class concerns the missing components in dynamical systems formulated as in \cite{harlim_2021}, where the underlying dynamics of the unresolved components is learned from the available data.
The second class arises in moment closure problems for kinetic equations.

Data-driven approaches to kinetic and moment closure have developed along several complementary directions, including learned hydrodynamic models for kinetic equations \cite{HanMaMaE2019}, gradient-based moment closures \cite{HuangChengChristliebRoberts2022}, structure- and invariance-preserving learned closures \cite{PorteousLaiuHauck2023,LiDongWang2023}, and hyperbolic learned closures for the BGK equation \cite{ChristliebDingHuangKrupansky2025}.
Related recent work has also investigated learning heat-flux closures directly from kinetic plasma data \cite{HuangDongWang2025,IngelstenEtAl2025}.
A recent overview of machine-learning approaches to plasma moment closure is provided in \cite{BurlesCamporeale2026}.

While machine-learning approaches, particularly neural-network-based models, have emerged as powerful tools for modeling dynamical systems, recent studies suggest that classical linear estimators, such as kernel ridge regression (KRR), can also produce accurate long-horizon predictions \cite{hamzi2021learning,darcy2025learning,song2025_KRR_DM,kreider2026learning}.
Despite differences in kernel choices and validation methodologies, a common theme in these studies is the importance of the validation procedure and the choice of validation metric used to determine the kernel hyperparameters.

One goal of this paper is to build upon these findings by investigating KRR for closure modeling.
Specifically, for the first class of closure problems described above, we replace the LSTM estimator used in \cite{harlim_2021} with a KRR estimator for the unresolved dynamical components.
Since the KRR solution can be computed efficiently when the training dataset is not excessively large, we investigate whether KRR can produce accurate long-horizon predictions using relatively small training datasets.
We compare its predictive performance with that of nonlinear neural-network-based estimators, such as LSTMs.

Because the modeling framework couples the learned closure model with a numerical approximation of the dynamical equations for the resolved components, we examine, both theoretically and numerically, the sensitivity of the prediction accuracy to the choice of numerical solver and to the manner in which the unresolved components are incorporated into the approximate solution operator for the resolved dynamics.

For the second class of problems, our goal is to reduce the computational cost of kinetic simulations.
Numerically resolving kinetic equations requires the expensive evolution of a distribution function in phase space, whereas the primary quantities of physical interest are often lower-dimensional macroscopic space-time fields, such as density, bulk velocity, temperature, and energy.
The evolution of these macroscopic fields can be described by a system of PDEs in physical space, thereby reducing the phase-space problem to a more manageable computational task.
However, the resulting PDE system is not closed because it depends on higher-order quantities, such as the stress tensor and heat flux.

Although standard treatments of kinetic and moment descriptions are available \cite{CercignaniIllnerPulvirenti1994,Levermore1996}, discrepancies between numerical fluxes constructed at the kinetic and macroscopic levels can introduce non-negligible errors, even when the underlying heat flux is fully specified.
These discrepancies may also depend on the choice of macroscopic solver.
Motivated by this observation, we propose using KRR to model the discrepancy between the flux generated by the kinetic solver and its approximation in the macroscopic solver.
Specifically, we learn the flux discrepancy as a function of the macroscopic variables at the same time instant, resulting in an algebraic closure model.
We evaluate the proposed closure approach on a family of smooth periodic thermally perturbed initial conditions, with the goal of predicting spatiotemporal fields arising from nontrivial multi-modal initial conditions.

The remainder of this paper is organized as follows.
In Section~\ref{sec:krr}, we provide an overview of learning missing components of dynamical systems using KRR and discuss the corresponding closure model and error analysis.
In Sections~\ref{sec:L63} and \ref{sec:KS}, we report numerical results for the Lorenz-63 and Kuramoto--Sivashinsky equations, respectively.
In Section~\ref{sec:kinetic}, we introduce the kinetic closure problem and present the numerical results.
We conclude the paper with a summary and a discussion of open problems in Section~\ref{sec:conclusion}.

Four appendices are included.
Appendix~\ref{app:model_selection} reviews the detailed validation strategy used to select the KRR hyperparameters.
Appendix~\ref{app:proof_of_theorem} provides the detailed proof of the theoretical results.
Appendix~\ref{app:KS_results} presents additional results for the Kuramoto--Sivashinsky example.
Finally, Appendix~\ref{app:kinetic} describes the training dataset used in the kinetic problems and reports additional numerical results.

\section{Learning missing components of dynamical systems with kernel ridge regression}
\label{sec:krr}

Consider dynamical systems for states $\bfx_n \in \bbR^{d_x}, \bfy_n \in \bbR^{d_y}$ that evolve based on the following maps:
    \begin{align}
\bfx_{n+1}
  &=
    \Phi_{\bfx} (\bfx_n,\bfy_n)
,
    \label{eq:x}
    \\
\bfy_{n+1}
  &=
    \Phi_{\bfy} (\bfx_n,\bfy_n)
,
    \label{eq:y}
    \end{align}
where $\Phi_{\bfx}$ and $\Phi_{\bfy}$ denote the solution operators of each component.
In this paper, we consider $\Phi_{\bfx}, \Phi_{\bfy}$ induced by the discretization of systems of ordinary and partial differential equations.
Here, $\bfx_n$ denotes the resolved components whereas $\bfy_n$ denotes the unresolved components.
With the structure in \eqref{eq:x}, we assume that we are either (1) given the time series of $\{\bfx_n,\bfy_n\}_{n=1}^N$; or (2) we can identify $\{\bfy_n\}_{n=1}^N$ from the observed $\{\bfx_n\}_{n=1}^N$.
Furthermore, $\Phi_{\bfx}$ is available in the following sense.
Since we are interested in identifying solutions of differential equations, the map $\Phi_{\bfx}$ denotes the approximate solution obtained by a numerical solver applied on a known system of ordinary differential equations, $\bfx' = f(\bfx,\bfy)$ resolved at $t_n=n\Delta t$.

Our task is to approximate the dynamical map $\Phi_{\bfy}$ of the unresolved components that will be assumed to be unknown (or missing) with a data-driven model.
We should point out that in many applications, the operator $\Phi_{\bfy}$ can sometimes depend on time-lagged data of $\bfx_n, \bfy_n, \ldots, \bfx_{n-m}, \bfy_{n-m}$ induced by some standard model reduction formalism \cite{harlim_2021}.
For simplicity of the discussion throughout this manuscript, we will discuss the case $m=0$ which corresponds to the Markovian dynamics in \eqref{eq:x}-\eqref{eq:y}.
We will demonstrate our findings for $m>0$ in the relevant numerical sections below.

In this paper, we will consider modeling \eqref{eq:y} using the recently introduced Diffusion Maps Kernel Ridge Regression \cite{song2025_KRR_DM}, which we will discuss in the following Subsection~\ref{sec:DMKRR}.
In practice, we should emphasize that we usually do not know the numerical solver that is used to generate the time series of $\bfx_n$ and $\bfy_n$.
This motivates us to study the errors of the predictions under the discrepancy of solvers being used to generate the data and to approximate the dynamics in Subsection~\ref{sec:closure_error}.

\subsection{Diffusion Maps Kernel Ridge Regression (DMKRR)}
\label{sec:DMKRR}

Given a training data set $\{\bfx_n \in \bbR^{d_x},\bfy_n\in \bbR^{d_y}\}_{{n=1}}^N$, define $\bfz_n = (\bfx_n^\top,\bfy_n^\top)^\top \in \bbR^{d_z}$, where $d_z=d_x+d_y$.
For any scalar-valued positive-definite kernel
    \[
k_{\varepsilon}
  :
    \bbR^{d_z}
    \times
    \bbR^{d_z}
  \rightarrow
    \bbR
,
    \]
with the kernel-scale parameter $\varepsilon$, we define the Gram matrix $K_\varepsilon \in \bbR^{N\times N}$ as follows:
    \[
[K_\varepsilon]_{i,j}
  =
    \left[
        k_\varepsilon
        \left(
            \bfz_i,\bfz_j
        \right)
    \right]
.
    \]
Kernel ridge regression approximates the $j^{\textnormal{th}}$ component of $(\Phi_{\bfy})_j(\bfz) \approx\sum_{i=1}^N k_\varepsilon(\bfz,\bfz_i)A_{i,j}$, where the coefficient matrix $A = [A_1,\ldots,A_{d_y}]\in\bbR^{N\times d_y}$ can be determined by solving
    \[
(A_{\varepsilon,\textnormal{reg}})_j
  =
    \arg\min_{A_{j}}
      \|Y_j-K_\varepsilon A_j\|^2_2
      +
      \lambda_{\textnormal{reg}}A_j^\top K_\varepsilon A_j
,\qquad
    \forall\; j= 1,\ldots, d_y
,
    \]
where $\|\cdot\|_2$ denotes the Euclidean norm, and $Y_j$ denotes the $j^{\textnormal{th}}$ column of $Y = [\bfy_1,\ldots, \bfy_N]^\top \in \bbR^{N \times d_y}$.
Here, $\lambda_{\textnormal{reg}}>0$ is the ridge regularization parameter.

In matrix form, the solution $A_{\varepsilon,\lambda_{\textnormal{reg}}} = [ (A_{\varepsilon,\lambda_{\textnormal{reg}}})_1 , \ldots , (A_{\varepsilon,\lambda_{\textnormal{reg}}})_{d_y} ]$ is given as,
    \[
A_{\varepsilon,\lambda_{\textnormal{reg}}}
  =
    \left(
        K_\varepsilon
        +
        \lambda_{\textnormal{reg}} I
    \right)^{-1}
    Y
.
    \]
where $I$ is the $N\times N$ identity matrix.

For a new input $\bfz$ and training data $Z = \{\bfz_1,\ldots, \bfz_N\}$, define
    \[
k_\varepsilon(Z,\bfz)
  =
    \begin{bmatrix}
        k_\varepsilon(\bfz_1,\bfz)
        & \cdots &
        k_\varepsilon(\bfz_N,\bfz)
    \end{bmatrix}^\top
    \in\bbR^{N}
.
    \]
The resulting KRR predictor is
    \[
\Phi_{\bfy,\varepsilon,\lambda_{\textnormal{reg}}}(\bfz)
  \coloneqq
    A_{\varepsilon,\lambda_{\textnormal{reg}}}^\top
    k_\varepsilon(Z,\bfz) \in \bbR^{d_y}
.
    \]
We note that each KRR model depends critically on the choice of kernel $k_\varepsilon$, bandwidth parameter $\varepsilon$, and the regularization parameter $\lambda_{\textnormal{reg}}$.
We will follow the same validation procedure proposed in \cite{song2025_KRR_DM} for choosing the hyperparameters $\{\varepsilon, \lambda_{\textnormal{reg}}\}$.
In particular, we will choose the parameters by fitting the predicted solutions obtained by rolling out (iterating) the KRR predictor to a validation trajectory (see Appendix~\ref{app:model_selection} for the details and the metric being used for validation).
We will denote the corresponding KRR predictor for optimal $\{\varepsilon^*, \lambda_{\textnormal{reg}}^*\}$ as,
    \begin{equation}
\widetilde{\Phi}_y
  \coloneqq
    \Phi_{y,\varepsilon^*,\lambda_{\textnormal{reg}}^*}
.
    \label{eq:KRRPredictor}
    \end{equation}

While one can use any kernel for the KRR fitting, we found the Diffusion Maps kernel \cite{song2025_KRR_DM,harlim2026diffusion} to give improved prediction over standard kernels such as the Gaussian kernel.

The standard radial basis function (RBF) kernel is
    \begin{equation}
\widetilde{k}_{\varepsilon}(\bfz,\bfz')
  =
    \exp
    \left(
        -\,\frac{\|\bfz-\bfz'\|_2^2}{4\varepsilon}
    \right)
,\qquad
    \varepsilon>0
.
    \label{eq:RBFKernel}
    \end{equation}

The Diffusion Maps kernel is a data-dependent kernel (referred to as DM) obtained by applying the following two normalization steps:
    \begin{alignat}{2}
q_{\varepsilon,N}(\bfz)
  &=
    \frac{1}{N}
    \sum_{j=1}^{N}
      \widetilde{k}_\varepsilon(\bfz,\bfz_j)
,\qquad
&
\widehat{k}_{\varepsilon,N}(\bfz,\bfz')
  &=
    \frac{
        \widetilde{k}_\varepsilon(\bfz,\bfz')
    }{
        q_{\varepsilon,N}(\bfz)
        q_{\varepsilon,N}(\bfz')
    }
,
  \notag
    \\
\widehat{q}_{\varepsilon,N}(\bfz)
  &=
    \frac{1}{N}
    \sum_{j=1}^{N}
      \widehat{k}_{\varepsilon,N}(\bfz,\bfz_j)
,\qquad
&
k_{\varepsilon,N}(\bfz,\bfz')
  &=
    \frac{
        \widehat{k}_{\varepsilon,N}(\bfz,\bfz')
    }{
        \sqrt{
            \widehat{q}_{\varepsilon,N}(\bfz)
            \widehat{q}_{\varepsilon,N}(\bfz')
        }
    }
.
    \label{eq:DMKernel}
    \end{alignat}
We refer to $k_{\varepsilon,N}$ in \eqref{eq:DMKernel} as the diffusion maps kernel.
One can show via similarity transform that the eigenvalues of $\frac{1}{N}K_{\varepsilon,N}$, where $K_{\varepsilon,N}$ denotes the Gram matrix of $k_{\varepsilon,N}$, are exactly identical to the eigenvalues of the resulting Diffusion Maps algorithm \cite{harlim2026diffusion}.
Under the manifold assumption, $\bfz\in\mathcal{M}$, the DM kernel converges to the heat kernel corresponding to the heat equation defined on the manifold $\mathcal{M}$ \cite{harlim2026diffusion}.

In our numerical study in Sections~\ref{sec:L63} and \ref{sec:KS}, we will compare the results using the Gaussian kernel in \eqref{eq:RBFKernel} and the DM kernel in \eqref{eq:DMKernel}.

\subsection{Closure model and error analysis}
\label{sec:closure_error}

Suppose that $\bfx(t)$, $\bfy(t)$ are the solution to a system of ODEs:
    \begin{equation}
\bfx'
  =
    \bff(\bfx,\bfy)
, \qquad
\bfy'
  =
    \bfg(\bfx,\bfy)
,
    \label{eq:ode}
    \end{equation}
under an initial condition $\bfx_0$, $\bfy_0$, where $f$ is given but $g$ is unknown, and they are Lipschitz continuous.
In practice, the time series data $\{\bfx_n,\bfy_n\}_{n=1}^N$ is obtained from an unknown numerical solver on $t_n=n\Delta t$.
In the following analysis, we will assume that the underlying solver is accurate such that $(\bfx_n=\bfx(t_n), \bfy_n=\bfy(t_n))_{n=1}^N$, and $t_n = n\Delta t$.

The proposed model is a discrete map,
    \begin{equation}
\widehat{\bfx}_{n+1}
  =
    \Phi_{\bfx} (\widehat{\bfx}_n,\widehat{\bfy})
,\qquad
\widehat{\bfy}
  =
    \mathcal{I}_n(\widehat{\bfy}_{n-p+1},\ldots,\widehat{\bfy}_{n+1})
,\qquad
\widehat{\bfy}_{n+1}
  =
    \widetilde{\Phi}_{\bfy} (\widehat{\bfx}_n,\widehat{\bfy}_n)
,
    \label{eq:estimator}
    \end{equation}
for $\widehat{\bfx}_0=\bfx_0$, $\widehat{\bfy}_0=\bfy_0$.
Here, $\Phi_{\bfx}$ denotes a numerical discretization of the first component in \eqref{eq:ode}.
Depending on the solver, which may be a multi-stage method such as RK4 or ETDRK4, one may need to evaluate the unresolved state $\widehat{\bfy}$ at intermediate times $t\in [t_n,t_{n+1}]$.
Since the learned map $\widetilde{\Phi}_y$ produces nodal values only at the discrete times $t_n=n\Delta t$, we introduce an interpolation operator $\mathcal{I}_n$ that uses the $p+1$ available nodal values $\widehat{\bfy}_{n-p+1}, \ldots, \widehat{\bfy}_{n+1}$ to approximate $\widehat{\bfy}(t)$ for $t \in [t_n,t_{n+1}]$.
In practice, we use a low-order Lagrange interpolant which is stable.

In the analysis below, we will study the effects of numerical discretization error in $\Phi_{\bfx}$ and the interpolation error associated with the stable interpolation operators $\mathcal{I}_n$, under the assumption that the KRR approximation error satisfies,
    \begin{equation}
\max_{0\leq n\leq N_T}
\| (\widetilde{\Phi}_y)^n(z_0) - (\Phi_{\bfy})^n(z_0) \|
  =
    \delta(\varepsilon,N) \to 0
\qquad\text{as}\qquad
    \varepsilon \to 0,
    \, N \to \infty
.
    \label{approxerror}
    \end{equation}
While the error for $n=1$ has been reported in \cite{harlim2026diffusion}, it is still an open problem to determine how $\delta$ scales as a function of $\epsilon$ and $N$, accounting for the validation strategy for choosing the hyper-parameters in the kernel and non i.i.d. data.
Under this assumption, one can deduce the following:

    \begin{theorem}
    \label{thm1}

Let the time series data $\{\bfx_n,\bfy_n\}_{n=1}^N$ be samples of the solution to \eqref{eq:ode} at $t_n = n\Delta t$, for some fixed $\Delta t>0$.
Let $\Phi_{\bfx}$ denote the one-step map induced by a globally order-$q$ numerical time stepping method for the $\bfx$-component of \eqref{eq:ode}.
Assume that the KRR predictor $\widetilde{\Phi}_y$, as defined in \eqref{eq:KRRPredictor}, satisfies \eqref{approxerror}.
Define the degree-$p$ stable interpolation operator, $\mathcal{I}_n: \bbR^{(p+1)\times d_y}\to C\big([t_n,t_{n+1}];\bbR^{d_y}\big)$.
For fixed $T = N_T \Delta t >0$,
    \[
\max_{0\leq n \leq N_T}
\|\bfx_n - \widehat{\bfx}_n\|
  \leq
    C_T
    [\delta(\varepsilon,N)+ (\Delta t)^{p+1} + (\Delta t)^q]
,
    \]
where $\widehat{\bfx}_n$ is defined by \eqref{eq:estimator} with initial conditions $\widehat{\bfx}_0=\bfx_0, \widehat{\bfy}_0=\bfy_0$.

    \end{theorem}

See the proof in Appendix~\ref{app:proof_of_theorem}.
In this case, one can see that, in the limit of $N \to \infty$, the dominant error depends on the order of the interpolation scheme and the ODE solvers.

\section{Lorenz--63 benchmark: An ODE example}
\label{sec:L63}

We consider the Lorenz--63 equations
    \begin{equation}
\dot{x}
  =
    \sigma(y-x)
,\qquad
\dot{y}
  =
    x(\rho-z)-y
,\qquad
\dot{z}
  =
    xy-\beta z
,
    \label{eq:L63}
    \end{equation}
with the standard parameters $\sigma =10$, $\rho=28$, $\beta =8/3$, for which the system exhibits chaotic dynamics.
Relating to the notation in previous section, we consider $\bfx = (x,y)$ as the resolved variable and $\bfy = z$ as the unresolved variable.
Here, $\Phi_{\bfx}$ denotes the map induced by approximating the first two components in \eqref{eq:L63} with the RK4 scheme with $\Delta t=0.01$.

Since we use RK4, we need to resolve $\widehat{\bfy}$ at intermediate time $t_n + \frac{\Delta t}{2}$.
Here, we employ degree-$p$ polynomial interpolation, where $p\in\{2,\ldots,5\}$.
As a lower-order reference, we also use a second-order Heun update in which the closure is evaluated only at the two endpoints, $\widehat{\bfy}_n$ and $\widehat{\bfy}_{n+1}$, so no interpolation is needed in this case.

\subsection{Experimental setup and evaluation metric}

We consider both the RBF and DM kernels defined in \eqref{eq:RBFKernel} and \eqref{eq:DMKernel} of Section~\ref{sec:DMKRR}.
We will construct the Gram matrix with training data sets of sizes $N_{\textnormal{train}} \in \{512,1024,2048,4096\}$.
For each configuration, ten independently trained models are constructed, each using a different training data set of size $N_{\textnormal{train}}$.
The kernel bandwidth and ridge parameter are selected using rollout-aware validation over $4096$ randomly sampled hyperparameter pairs, with three validation rollouts used to select each candidate model.
Each selected model is subsequently evaluated on $500$ test initial conditions.

To evaluate the prediction skill, we consider a metric referred to as the Valid Prediction Time (VPT), used in \cite{mandal2025learning,song2025_KRR_DM,kreider2026learning}, defined as follows:
    \begin{equation}
\textnormal{VPT}
  =
    \Lambda_1
    \left(
        t_{n_*}-t_0
    \right)
,\qquad
n_*
  =
    \max
    \left\{
        n
        :
        \mathcal{E}_i\leq\gamma
        \text{ for all }i\leq n
    \right\}
,
    \label{eq:L63VPT}
    \end{equation}
where $\Lambda_1 = 0.91$ denotes the leading Lyapunov exponent (effectively, $\Lambda_1$ normalizes the VPT in Lyapunov-time units), and
    \begin{equation}
\mathcal{E}_i
  =
    \left[
        \frac{1}{d_z}
        \sum_{j=1}^{d_z}
          \left(
              \frac{
                  \widehat{\bfz}_{i,j}-\bfz_{i,j}
              }{
                  s_j
              }
          \right)^2
    \right]^{1/2}
.
    \label{eq:L63NormalizedError}
    \end{equation}
Here, $\widehat{\bfz}_{i,j}$ and $\bfz_{i,j}$ denote the $j^{\textnormal{th}}$ components of $\widehat{\bfz}=(\widehat{\bfx},\widehat{\bfy})$ and  $\bfz=(\bfx,\bfy)$ at time index $i$, respectively, and $s_j$ denotes the empirical standard deviation of the $j^{\textnormal{th}}$ component of $\bfz$.
Intuitively, VPT measures the length of time, in Lyapunov-time units, over which the normalized prediction error $\mathcal{E}_i$ remains below the threshold $\gamma$.
In this numerical experiment, we set $\gamma=0.3$.

For comparison, we also train a hybrid closure model based on an LSTM model, with sequence length $m=19$ and rollout length $R=5$, where $m$ denotes the number of most recent states used as input and $R$ denotes the number of autoregressive prediction steps applied during training.
The LSTM takes an input dimension of $d_z \times m$ (where $d_z = 3$ represents $x,y,z$) and a hidden state dimension of $h=500$, followed by a fully connected layer that maps the final hidden state to a scalar prediction.
The model is trained for $T_{\max} = 500$ epochs with a batch size of 256 using the Adam optimizer~\cite{kingma2014adam}.
The learning rate $\eta_t$ at epoch $t$ follows a cosine annealing schedule, decaying from an initial value of $\eta_{\max} = 10^{-3}$ down to a minimum of $\eta_{\min} = 10^{-6}$, calculated as: $\eta_t = \eta_{\min} + \frac{1}{2}(\eta_{\max} - \eta_{\min})\left(1 + \cos\left(\frac{t}{T_{\max}}\pi\right)\right)$.
During inference, each prediction step requires one LSTM forward pass over a length-$m$ sequence, together with one numerical update such as Cubic-RK4.
Specifically, at each of the $m$ time steps, the LSTM gates require multiplying the input and hidden state vectors by matrices of size $h \times d_z$ and $h \times h$, respectively.
Thus, the dominant per-step cost is $O(m(hd_z+h^2))$, which readily reduces to $O(mh^2)$ since $h \gg d_z$.

\subsection{Results}

Table~\ref{table:L63InterpolationOrder} compares the stage-reconstruction strategies at the largest training size, $N_{\textnormal{train}}=4096$.
The endpoint-only Heun coupling produces a VPT of approximately $4.2$ Lyapunov times for both kernels.
Quadratic ($p=2$) reconstruction increases the VPT to approximately $5.9$, while cubic reconstruction increases it further to approximately $8$ Lyapunov times.

Quartic ($p=4$) and quintic ($p=5$) interpolations provide modest additional gains beyond cubic ($p=3$) order, with the best mean value obtained by quintic DM-KRR.
These results indicate that while a higher-order ODE solver is critical for more accurate prediction, it is also important to suppress the prediction error with an interpolator that is as accurate as the ODE solver order, as suggested by Theorem~\ref{thm1}.

The VPT of the KRR models with $p\geq 3$ nearly doubles that of the corresponding LSTM baseline at $N_{\textnormal{train}}=4096$, while also exhibiting substantially less trained-model variability.
The difference between the RBF and diffusion kernels is comparatively small: both kernels benefit strongly from improved stage reconstruction, and the DM kernel yields smaller variability.

Table~\ref{table:L63CubicTrainingSize} reports results for cubic reconstruction across the four training-set sizes.
The KRR performance is relatively insensitive to $N_{\textnormal{train}}$ over the range considered.
Apart from increased variability at $N_{\textnormal{train}}=1024$, both kernels maintain mean VPT values near 8 Lyapunov times.
In contrast, the LSTM baseline remains near $4.2$-$4.6$ Lyapunov times across the training sizes.
For this example, improving the numerical coupling between the learned closure and the resolved integrator has a substantially larger effect than increasing the number of training samples within the tested range.

The relatively weak dependence on training-set size and kernel choice shows that, for this problem, the numerical consistency of the hybrid time integrator is at least as important as the supervised regression method itself, as suggested by Theorem~\ref{thm1}.

\begin{table}[t]
    \centering
    \begin{tabular}{lccc}
        \toprule
        Stage treatment
        &
        DM-KRR
        &
        RBF-KRR
        &
        LSTM
        \\
        \midrule
        Heun / endpoint RK2
        &
        $4.1969 \pm 0.0799$
        &
        $4.1783 \pm 0.1165$
        &
        $4.4603 \pm 1.3410 $
        \\
        Quadratic RK4
        &
        $5.9874 \pm 0.0868$
        &
        $5.9242 \pm 0.0980$
        &
        $4.4924\pm 1.2797$
        \\
        Cubic RK4
        &
        $8.3815 \pm 0.3321$
        &
        $7.9683 \pm 0.4844$
        &
        $4.6008 \pm 1.3046$
        \\
        Quartic RK4
        &
        $8.4832 \pm 0.4910$
        &
        $8.3477 \pm 0.6612$
        &
        $4.6115 \pm 1.3114$
        \\
        Quintic RK4
        &
        $\mathbf{8.6505 \pm 0.5239}$
        &
        $8.4533 \pm 0.7001$
        &
        $4.6439 \pm 1.2873$
        \\
        \bottomrule
    \end{tabular}
    \caption{
        Mean test valid prediction time for the Lorenz--63 hybrid closure models with $N_{\textnormal{train}}=4096$.
        The KRR hyperparameters are selected using three rollout validations per candidate pair.
        Results are reported as the mean and standard deviation across independently trained-model test means.
    }
    \label{table:L63InterpolationOrder}
\end{table}

\begin{table}[t]
    \centering
    \begin{tabular}{cccc}
        \toprule
        $N_{\textnormal{train}}$
        &
        DM-KRR
        &
        RBF-KRR
        &
        LSTM
        \\
        \midrule
        512
        &
        $8.1719 \pm 0.4849$
        &
        $7.8863 \pm 0.5426$
        &
        $4.2784 \pm 1.3443$
        \\
        1024
        &
        $7.4267 \pm 1.3808$
        &
        $7.4753 \pm 0.9587$
        &
        $4.2253 \pm 1.3061$
        \\
        2048
        &
        $8.2535 \pm 0.4977$
        &
        $8.0014 \pm 0.4216$
        &
        $4.2652 \pm 1.2222$
        \\
        4096
        &
        $\mathbf{8.3815 \pm 0.3321}$
        &
        $7.9683 \pm 0.4844$
        &
        $4.6008 \pm 1.3046$
        \\
        \bottomrule
    \end{tabular}
    \caption{
        Mean test valid prediction time for the Lorenz--63 hybrid models using cubic reconstruction of the closure at RK4 stage times.
        Results are reported as the mean and standard deviation across independently trained-model test means.
    }
    \label{table:L63CubicTrainingSize}
\end{table}

\section{Kuramoto--Sivashinsky benchmark: A PDE example}
\label{sec:KS}

We consider the one-dimensional Kuramoto--Sivashinsky equation on an $L$-periodic spatial domain, where $L=\frac{2\pi}{\sqrt{0.085}}$.
In Fourier coordinates, the PDE equation consists of an infinite-dimensional system of ODEs,
    \[
\frac{d v_k}{dt}
  =
    \left(
        q_k^2-q_k^4
    \right)v_k
    -
    \frac{iq_k}{2}
    \sum_{\ell\in\mathbb Z}
      v_\ell v_{k-\ell}
,
    \]
where $q_k=\frac{2\pi k}{L}$.
The linear term contains both the long-wave instability and short-wave dissipation, while the quadratic convolution transfers energy among Fourier modes.

Reference trajectories are generated using a Fourier pseudospectral discretization with dynamics truncated to $1\leq |k| \leq 48$ and $3/2$-dealiasing.
Time integration is performed with the fourth-order exponential time-differencing Runge--Kutta method (ETDRK4) with an integration time step $\delta t=0.005$.
For our study, the training data is subsampled at every 10 steps, or $\Delta t=0.05$.
These full-order trajectories are treated as the reference dynamics from which both the resolved variables and unresolved closure terms are extracted.

\subsection{Resolved model and unresolved closure}

Following the setup in \cite{harlim_2021}, we define the leading six $d_x=6$ Fourier modes $\bfx= (v_1,\ldots,v_{d_x})^\top\in\bbC^{d_x}$ as the resolved components.
For $k=1,\ldots,d_x$, define the resolved part of the quadratic nonlinearity by
    \[
B_k(\bfx)
  =
    -
    \frac{iq_k}{2}
    \sum_{
        \substack{
            1\leq |\ell|\leq d_x
            \\
            1\leq |k-\ell|\leq d_x
        }
    }
      v_\ell v_{k-\ell}
.
    \]
The resolved dynamics can then be written as
    \begin{equation}
\frac{dv_k}{dt}
  =
    \left(
        q_k^2-q_k^4
    \right)v_k
    +
    B_k(\bfx)
    +
    \theta_k
,\qquad
    k=1,\ldots,d_x
,
    \label{eq:KSResolvedDynamics}
    \end{equation}
where
    \[
\theta_k
  =
    -
    \frac{iq_k}{2}
    \sum_{\ell\in\mathbb Z}
      v_\ell v_{k-\ell}
    -
    B_k(\bfx)
    \]
is the contribution of all nonlinear interactions not contained in the resolved convolution.
In particular, $\theta_k$ includes interactions in which at least one participating mode lies outside the retained range.
The map $\Phi_{\bf x}$ is defined as the solution operator to \eqref{eq:KSResolvedDynamics} induced by ETDRK4 scheme with $\delta t=0.005$.
We define $\bfy = (\theta_1,\ldots,\theta_{d_y})^\top\in \bbC^{d_y}$ as the unresolved component, where $d_y=d_x=6$.
In principle, the dynamical evolution for $\bfy$ can be derived using the Mori-Zwanzig formalism \cite{harlim_2021}, resulting in an integro-differential equation whose discrete approximation yields a map
    \[
(\bfz_{n-m},\ldots,\bfz_n)
  \overset{\Phi_{\bfy}}{\longmapsto}
    \bfy_{n+1}
,
    \]
where $\bfz = (\bfx,\bfy)$ for some memory length $m>0$.
Here $t_n=n\Delta t$.
Since $\Phi_{\bfx}$ takes $\widehat{\bfz} =(\widehat{\bfx},\widehat{\bfy})$ solutions at higher resolution $\delta t = \Delta t/10$, we employ a linear interpolator $\mathcal{I}$ to attain estimates of $\widehat{\bfy}$ at $[t_n, t_{n+1}]$.

In our numerical experiment below, we shall see that with the KRR approximation $\widetilde{\Phi}_{\bfy}$, we do not gain a consistent advantage with the non-Markovian approximation, $m>0$.
Indeed, this approximation is competitive with a LSTM approximation that uses a large memory $m=19$ and is trained on a much larger training data set.

\subsection{Experimental setup and baseline}

In our implementation, $\bfx$ and $\bfy$ are turned into vectors of real numbers by concatenating the real and imaginary components of each complex vector.

\begin{figure}[t]
    \centering
    \includegraphics[
        width=0.45\textwidth
    ]{
        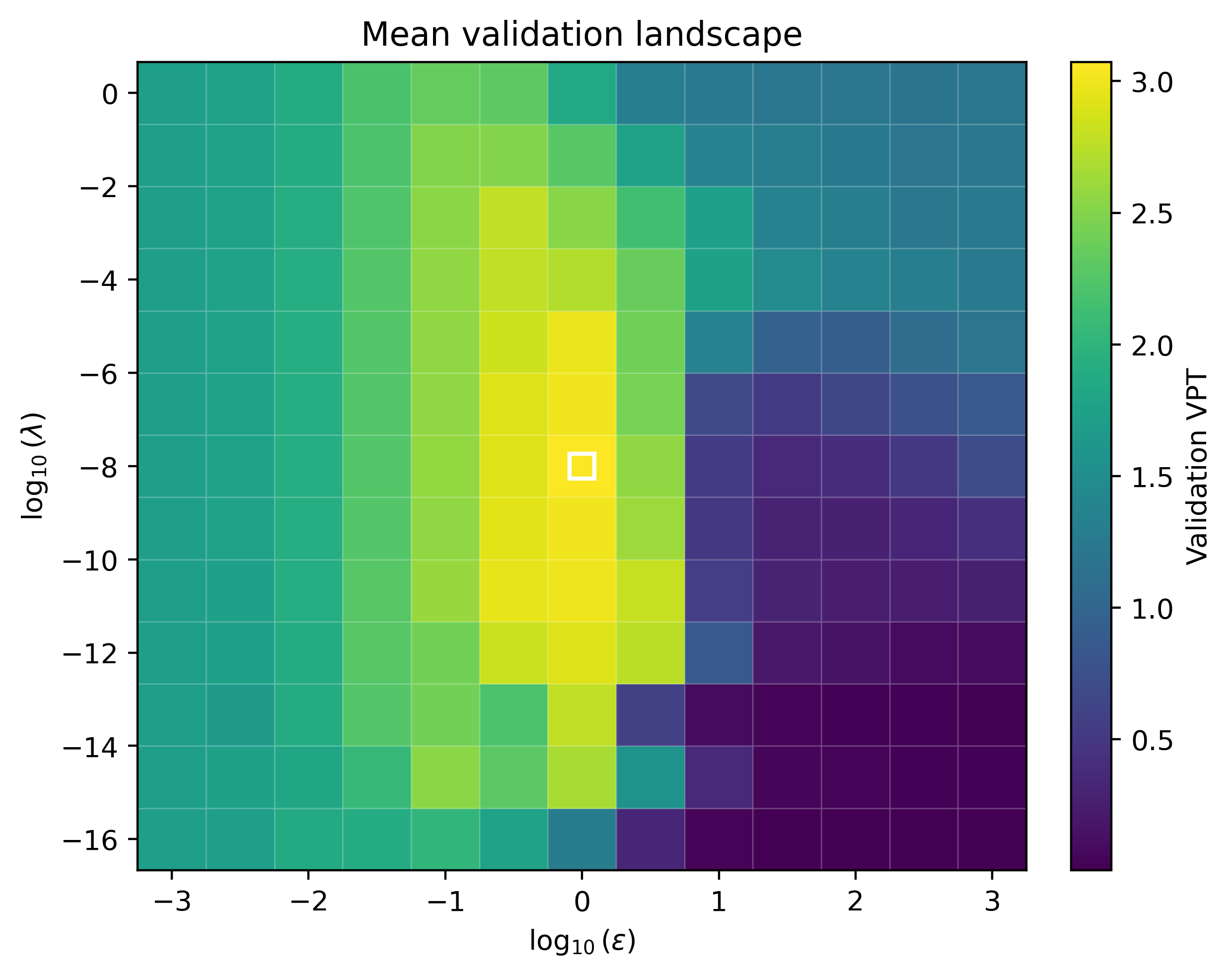
    }
    \includegraphics[
        width=0.45\textwidth
    ]{
        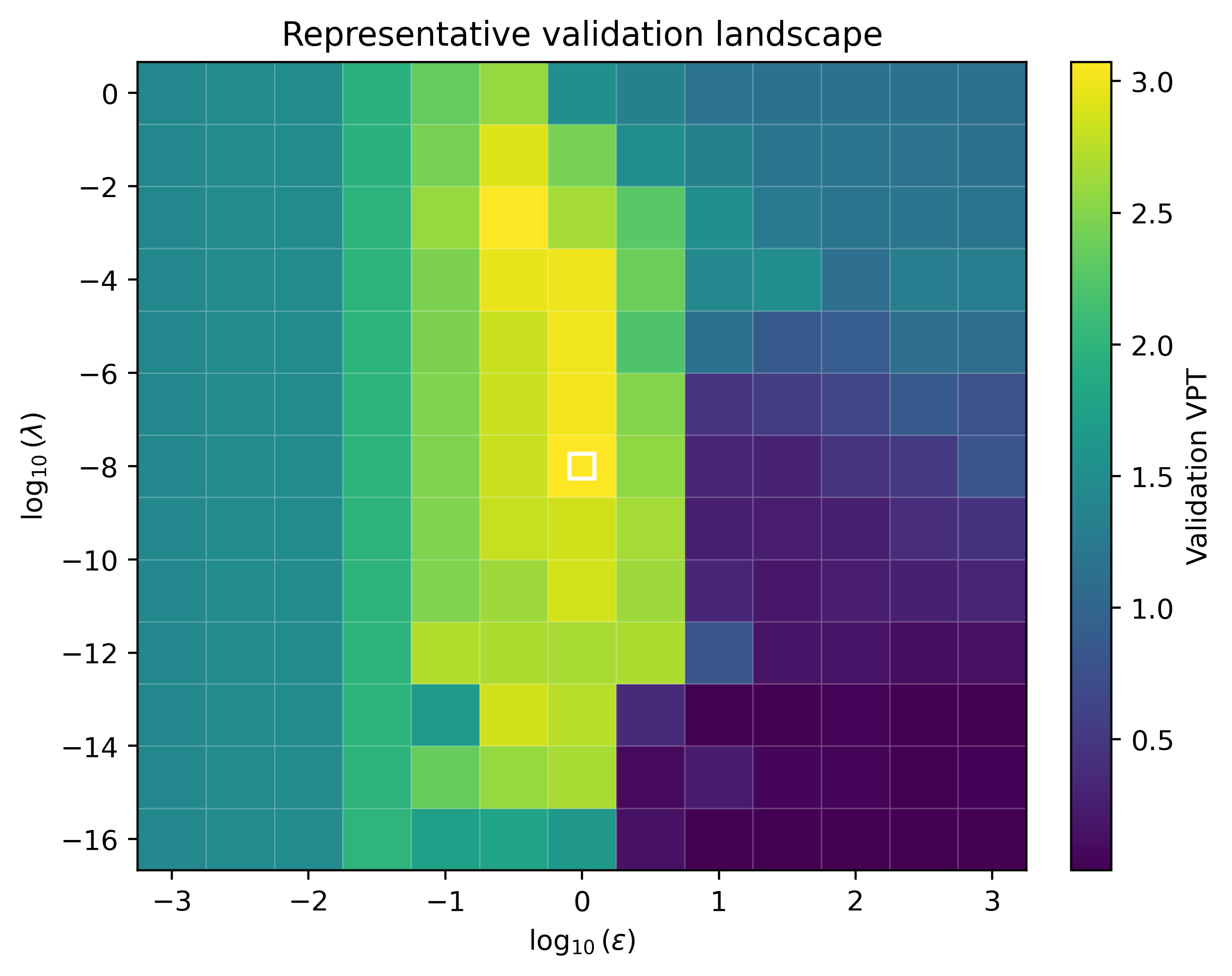
    }
    \caption{
        Mean validation VPT (left) and a representative validation (right) over the diffusion-kernel hyperparameter grid for $N_{\textnormal{train}}=16384$ and $m=0$, over ten training windows.
        The figure shows the mean landscape over ten independent training windows.
        The outlined square marker identifies the selected grid point.
    }
    \label{fig:KSValidationLandscape_mean}
\end{figure}

While the validation process follows the procedure described in Appendix~\ref{app:model_selection}, the bandwidth $\varepsilon$ and ridge parameters $\lambda_{\textnormal{reg}}$ for the kernels are selected from the deterministic logarithmic grid
    \[
\log_{10}\varepsilon
    \in[-3,3]
,\qquad
\log_{10}\lambda_{\textnormal{reg}}
    \in[-16,0]
,
    \]
using 13 values in each direction, for a total of $13\times 13=169$ validation cases per trained model.
The use of such a deterministic grid reduces the number of matrix inversions required in the training procedure.
In Figure~\ref{fig:KSValidationLandscape_mean}, we show the mean validation landscape and a representative validation landscape.
One can see that there is a smooth and local region where the VPT is high over the proposed parameter window.
We identify the selected parameter with a square marker.

For each training size $N_{\textnormal{train}}\in\{2048,4096,8192,16384\}$, ten independently trained KRR models are constructed from disjoint training windows.
Hyperparameters are selected using a separate validation window of length $2N_V$, with $N_V=2000$.
Five validation rollouts of length $N_V$, distributed across this window, are used for each candidate hyperparameter pair.
Each selected model is then evaluated on 50 recursive test rollouts.

The VPT metric as defined in \eqref{eq:L63VPT} is used to validate and measure the prediction skill.
Here, the normalized error metric in \eqref{eq:L63NormalizedError} is averaged over the resolved components, $\bfx$, rather than $\bfz$.
For the numerical experiment in this section, we set the threshold parameter in the \eqref{eq:L63VPT} to be $\gamma=0.5$.
The leading Lyapunov exponent is $\Lambda_1\approx 0.04$.

\paragraph{Baseline models.}

The no-closure baseline model sets $\bfy = 0$ and yields a VPT of $1.0815\pm0.4570$.
This will serve as a minimum baseline VPT to measure the prediction skill against.
On the other hand, the oracle closure model, which corresponds to the model that takes $\widehat{\bfy}$ attained from interpolating the true trajectory $\{\bfy_n\}$, has VPT $9.2943\pm 1.0185$.
The oracle model error is only induced by the interpolation step in \eqref{eq:estimator} and its VPT serves as an upper bound of the prediction skill.

We adopted a customized RNN architecture from prior work~\cite{harlim_2021} for sequence prediction.
The model uses an LSTM model with $24\times (m+1)$-dimensional inputs, a $100$-dimensional hidden state, and $12$-dimensional outputs.
Instead of standard linear gating, it employs three-layer MLPs with ReLU activations for the input-to-hidden and hidden-to-hidden transformations, mapping representations to a 1000-dimensional space before computing the four LSTM gates.
At each time step, the hidden state is passed through a two-layer fully connected neural network decoder that reduces it to $50$ dimensions before producing the final $12$-dimensional prediction.
The VPTs for various sizes of training data sets and memory $m$ are reported in Table~\ref{table:lstm_etdrk4} in Appendix~\ref{app:KS_results}.

\subsection{Results}

\begin{figure}[htbp]
    \centering
    \includegraphics[
        width=0.45\textwidth
    ]{
        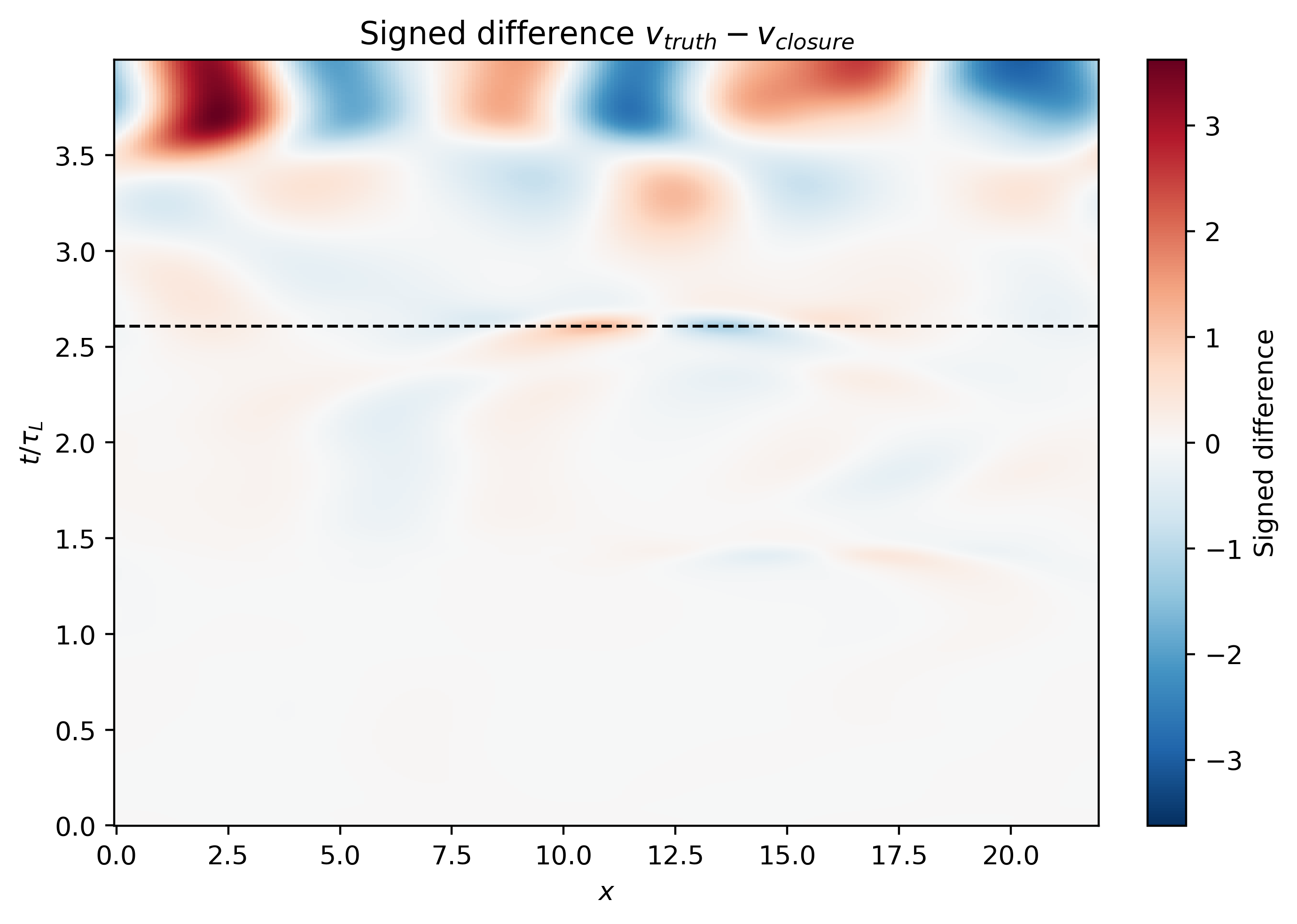
    }
    \includegraphics[
        width=0.45\textwidth
    ]{
        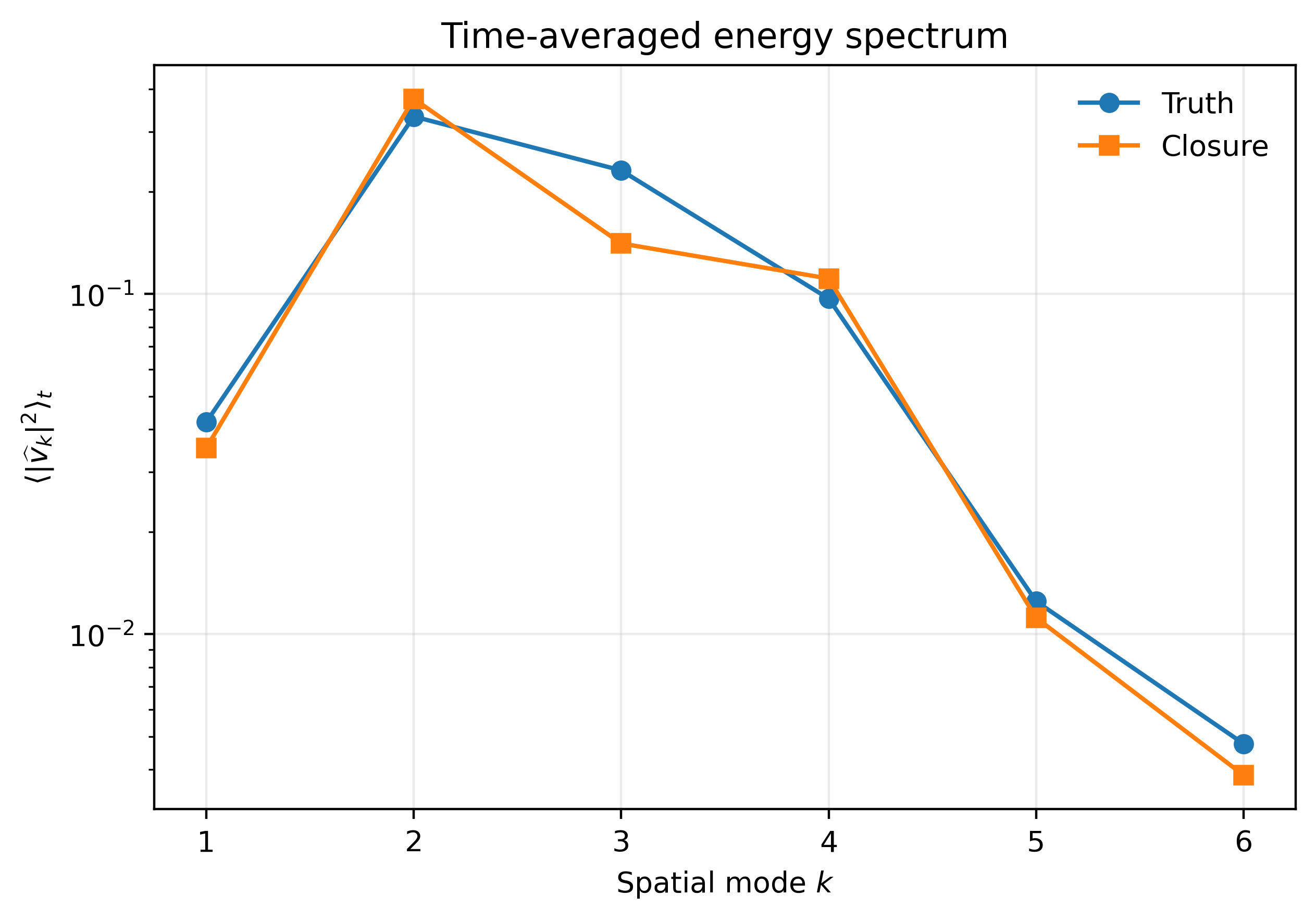
    }
    \caption{
        Left: Representative test rollout for the
        $N_{\textnormal{train}}=16384$, $m=0$, diffusion-kernel model.
        The signed physical-space difference between the full-order resolved field and the hybrid closure rollout is recorded.
        The dashed line marks the verified VPT,
        $2.608$.
        Right: The time-averaged energy in the six retained Fourier modes is shown.
    }
    \label{fig:KSRepresentativeRollout}
\end{figure}

Figure~\ref{fig:KSRepresentativeRollout} presents a representative test trajectory for the diffusion maps kernel model with $N_{\textnormal{train}}=16384$ and $m=0$.
The left panel shows the difference between the reference resolved field and the hybrid closure rollout.
The field is reconstructed on a dense spatial grid using the six retained Fourier modes.
The dashed line marks the verified VPT, $2.608$, which is longer than the VPT of the benchmark LSTM with $m=19$ and trained on much larger data set $2.5\times 10^5$ as reported in \cite{harlim_2021}.
Before this time, the trajectory difference remains comparatively small and spatially coherent.
After the threshold crossing, the difference grows rapidly, as expected for a chaotic system.

The right panel compares the time-averaged energy of the six resolved Fourier modes.
The learned closure reproduces the overall modal-energy distribution and the decay toward the higher retained modes.
The largest visible discrepancy occurs in the intermediate resolved modes, particularly near $k=3$, while the lowest and highest retained modes remain comparatively close.
Thus, although trajectory-level phase agreement is necessarily finite, the learned closure also captures important aspects of the resolved spectral structure over the displayed interval.

\begin{figure}[htbp]
    \centering
    \includegraphics[
        width=0.6\textwidth
    ]{
        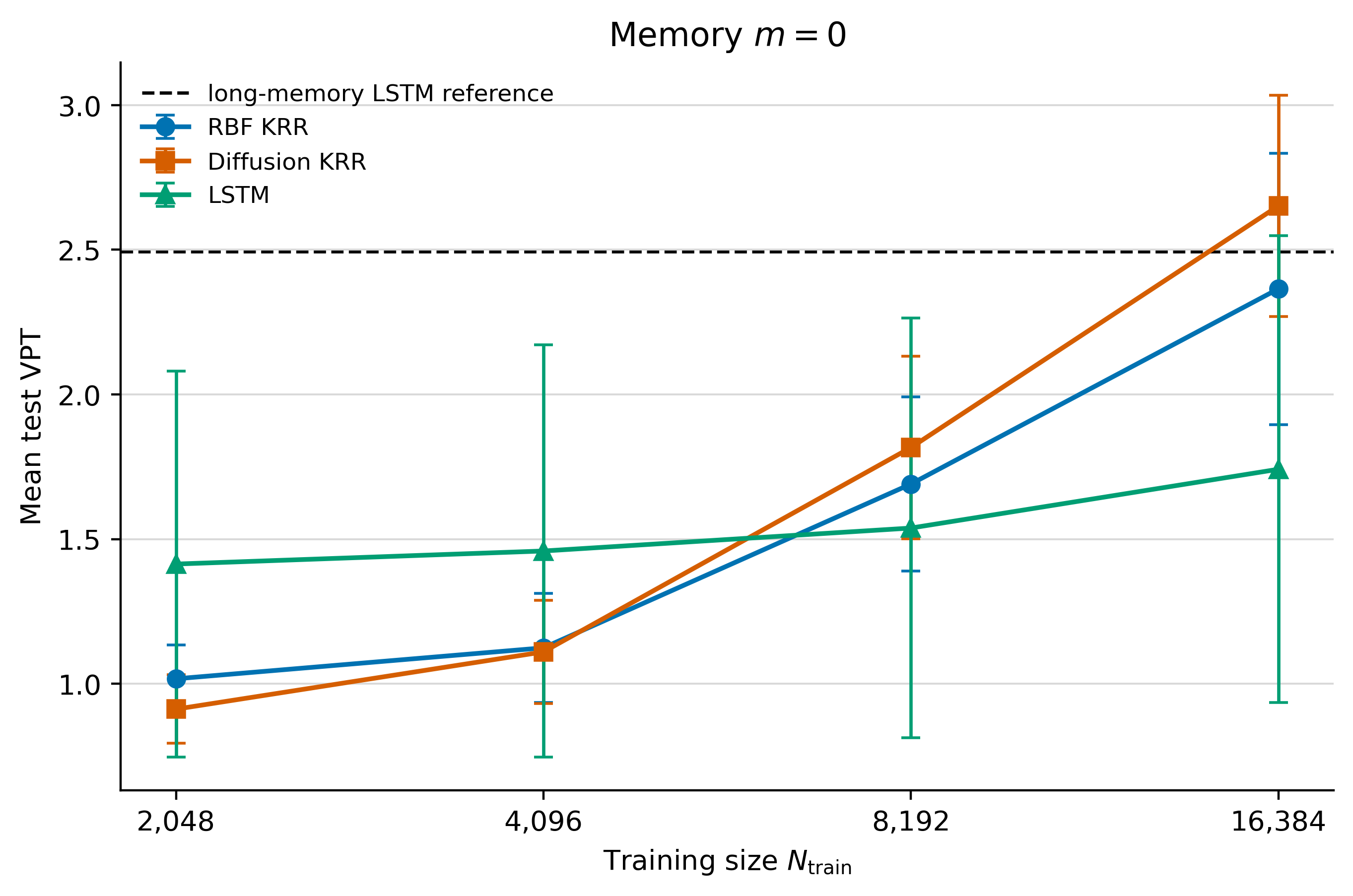
    }
    \caption{
        Mean test VPT as a function of training-set size for the $m=0$ KS closure models.
        KRR error bars denote variability across independently trained models.
        The LSTM values are included as a qualitative comparison because the LSTM training-data convention and noise augmentation are not identical to those used for KRR.
        The dashed reference line marks the reported long-memory LSTM mean, $2.494688$, whose reported standard deviation is $0.95904$.
    }
    \label{fig:KSTrainingSize}
\end{figure}

In Figure~\ref{fig:KSTrainingSize}, we report more comprehensive comparison statistics.
These statistics are obtained from ten independently trained-model test means, with 50 test rollouts per model.
Overall, these statistics are computed over 500 samples.
In addition to the baseline LSTM model trained in $m=19$ and a much larger data set, $2.5\times 10^5$, we also report the prediction of the LSTM model using the same size of training data.
Particularly, since each kernel predictor uses an additional $2N_V = 4000$ data for validation, in addition to using $N_{\text{train}}$ for constructing the Gram matrix, each LSTM prediction shown in Figure~\ref{fig:KSTrainingSize} corresponds to training on datasets of size $N_{\text{train}}+2N_V$.
We also report the corresponding quantitative results in Tables~\ref{table:KSTrainingSizeMemory} and \ref{table:lstm_etdrk4} in Appendix~\ref{app:KS_results}.

We note that the prediction skills of both kernel models are not better than that of the no-closure baseline ($1.0815\pm0.4570$) when $N_{\text{train}} \leq 4096$.
They are significantly better than this baseline when the training data is at least $N_{\text{train}}=8192$, and are competitive to the LSTM reference trained on $m=19$ and much larger data set, $2.5\times 10^5$, when $N_{\text{train}}=16384$.
In this regime, the DM kernel provides a slightly higher VPT than the Gaussian RBF kernel.
The LSTM models trained on a similar size dataset to the kernel model produce better predictions on small datasets but their skill increases slowly as a function of $N_{\text{train}}$; their best VPT, trained on $N_{\text{train}}+2N_V=16384+4000$, is only $1.741\pm 0.807$.

We should point out that for $2048\leq N_{\textnormal{train}}\leq16384$, we found that increasing from $m=0$ to $m=1$ (and even longer memory, $m>1$) does not improve prediction skill of both the kernel and LSTM models (see Tables~\ref{table:KSTrainingSizeMemory} and \ref{table:lstm_etdrk4}).

\section{Application to a kinetic closure problem}
\label{sec:kinetic}

We next apply the partially data-driven framework of Section~\ref{sec:DMKRR} to a kinetic-to-macroscopic closure problem.
In particular, we are interested in the closure of dynamical equations for the lower-dimensional macroscopic space-time fields, including density, bulk velocity, temperature, and energy.
Unlike the unresolved dynamical variables learned in the L63 and KS examples of Sections~\ref{sec:L63} and \ref{sec:KS}, the effective flux correction will be identified as a function of the current macroscopic state at each step, resulting in an algebraic closure for the macroscopic solver.

The remainder of this section is organized as follows.
In Section~\ref{sec:kinetic_setup}, we present the problem formulation.
Then, we discuss the proposed closure models and the training configuration in Sections~\ref{sec:kinetic_representations} and \ref{sec:kinetic_evaluation}, respectively.
In Section~\ref{sec:kinetic_results}, we present numerical results.

\subsection{Kinetic benchmark and solver-consistent closure}
\label{sec:kinetic_setup}

We consider the one-dimensional BGK equation, which replaces the expensive Boltzmann collision operator by relaxation toward a local Maxwellian while preserving the collision invariants \cite{BhatnagarGrossKrook1954,CercignaniIllnerPulvirenti1994}:
    \begin{equation}
\partial_t f + v \partial_x f
  =
    \frac{1}{\textnormal{Kn}}
    \left(
        M[f]-f
    \right)
,
    \label{eq:bgk_1d}
    \end{equation}
where $\textnormal{Kn}$ is the Knudsen number, $f=f(x,v,t)$ is a unit-mass particle distribution, and $M[f]$ is the local Maxwellian having the same first three velocity moments as $f$,
    \begin{equation}
M[f]
  =
    \frac{\rho}{\sqrt{2\pi\theta}}
    \exp
    \left(
        -\frac{(v-u)^2}{2\theta}
    \right)
.
    \label{eq:bgk_maxwellian}
    \end{equation}
Here, $\rho$, $u$, and $\theta$ denote the mass density, bulk velocity, and temperature, respectively; in the present nondimensionalization, $\theta$ is also the variance of the velocity distribution.
With $\psi(v) = (1,v,\frac{1}{2}v^2)^\top$, the conserved macroscopic state $\bfU$, energy $E$, and heat flux $q$ are
    \begin{equation}
\bfU
  =
    \int_{\bbR}\psi(v)f\,\diff v
  =
    \begin{pmatrix}
        \rho\\
        \rho u\\
        E
    \end{pmatrix}
,\qquad
E
  =
    \frac{1}{2}\rho u^2
    +
    \frac{1}{2}\rho\theta
,\qquad
q
  =
    \frac{1}{2}\int_{\bbR}(v-u)^3f\,\diff v
.
    \end{equation}
Taking moments of \eqref{eq:bgk_1d} gives
    \begin{equation}
\partial_t \bfU
+
\partial_x
\bfG(\bfU,q)
  =
    \mathbf{0}
,\qquad
\bfG(\bfU,q)
  =
    \begin{pmatrix}
        \rho u\\
        \rho u^2+\rho\theta\\
        (E+\rho\theta)u+q
    \end{pmatrix}
.
    \label{eq:kinetic_macro_system}
    \end{equation}
Thus, at the continuum level, the heat flux $q$ is the closure information not determined by the conserved state $\bfU$.
More generally, obtaining closed systems for finitely many velocity moments is the classical moment-closure problem; see, for example, \cite{Levermore1996}.

The benchmark consists of smooth periodic thermal perturbations on $(x,v)\in[-1,1]\times[-8,8]$.
The initial macroscopic state is prescribed by
    \begin{equation}
\rho_0(x)
  =
    1
,\qquad
u_0(x)
  =
    0
,\qquad
\theta_0(x)
  =
    1
    +
    \sum_{\ell=1}^{P}
      A_\ell
      \exp
      \left[
          -
          \left(
              \frac{d_{\textnormal{per}}(x,\mu_\ell)}{\sigma_\ell}
          \right)^2
      \right]
,\qquad
    P\in\{1,2\}
,
    \label{eq:macroscopic_ic_fields}
    \end{equation}
where
    \[
d_{\textnormal{per}}(x,\mu)
  =
    \min_{k\in\bbZ}
      |x-\mu+2k|
,\qquad
\mu_\ell
  =
    (j_\ell+\eta_\ell)\Delta x
,\qquad
j_\ell
  \in
    \bbZ
,\quad
\eta_\ell
  \in
    [-1/2,1/2)
.
    \]
The corresponding kinetic distribution is initialized using the conservative discrete Maxwellian described in Appendix~\ref{app:kinetic_solver}, whose discrete velocity moments reproduce the prescribed macroscopic state \eqref{eq:macroscopic_ic_fields} to the nonlinear solver tolerance.
The case $P=1$ defines the unimodal family, while $P=2$ defines the asymmetric bimodal family.
The amplitudes $A_\ell$, widths $\sigma_\ell$, phases $\mu_\ell$, and training/evaluation assignments are specified in Section~\ref{sec:kinetic_evaluation}.

Reference solutions of \eqref{eq:bgk_1d} are generated using a second-order finite-volume discretization and a second-order implicit-explicit (IMEX) Runge--Kutta (RK) time-stepping method \cite{PareschiRusso2005,habbershaw2025year4}.
Unless otherwise stated, we use $N_x=256$, $N_v=128$, $x\in[-1,1]$, $v\in[-8,8]$, and $\textnormal{Kn}=10^{-3}$.
The time step is specified using a hyperbolic CFL-type condition as $\Delta t = 4.39453125\times10^{-4}$, and each trajectory contains 456 time nodes, with final saved time $t_{\textnormal{F}}=0.199951171875\approx0.2$.

Here and below, $n$ denotes the time-step index, $s$ the Runge--Kutta stage index, $i$ the spatial cell index (so that $i+1/2$ denotes a spatial interface), and $k$ the discrete-velocity index.

At each Runge--Kutta stage, the kinetic solver first reconstructs the distribution function at each cell interface and applies a velocity-upwind method.
If $f_{i+1/2,k}^{-,n,s}$ and $f_{i+1/2,k}^{+,n,s}$ denote the reconstructed traces from the left and right at velocity $v_k$, respectively, then the kinetic distribution flux is
    \begin{equation}
\mathcal{F}_{i+1/2,k}^{n,s}
  =
    v_k^+
    f_{i+1/2,k}^{-,n,s}
    +
    v_k^-
    f_{i+1/2,k}^{+,n,s}
,\qquad
v_k^\pm
  =
    \frac{1}{2}
    \left(
        v_k\pm|v_k|
    \right)
.
    \label{eq:kinetic_distribution_flux}
    \end{equation}
Taking its conserved velocity moments gives the kinetic numerical interface flux
    \begin{equation}
\mathcal{G}_{\textnormal{kin},i+1/2}^{n,s}
  =
    \sum_k
      \omega_k
      \psi(v_k)
      \mathcal{F}_{i+1/2,k}^{n,s}\in\bbR^3
,
    \label{eq:kinetic_numerical_moment_flux}
    \end{equation}
where $\omega_k$ denotes the velocity quadrature weight.

On a finite velocity grid, evaluating the continuum Maxwellian \eqref{eq:bgk_maxwellian} at the discrete velocities does not, in general, preserve the prescribed mass, momentum, and energy under the discrete velocity quadrature.
To overcome this issue, we embed the conservative discrete Maxwellian of Mieussens~\cite{Mieussens2000} in the kinetic solver.
Its parameters are determined by solving a nonlinear algebraic system so that its discrete mass, momentum, and energy moments agree with those of the kinetic distribution to nonlinear-solver tolerance.
As a result, the discrete BGK collision term has vanishing conserved moments to this tolerance.
Further details are given in Appendix~\ref{app:kinetic_solver}.

The macroscopic solver instead constructs an interface flux directly from the reconstructed macroscopic states.
Its zero-heat-flux Rusanov (local Lax--Friedrichs) numerical flux \cite{Rusanov1962,Toro2009} has the form
    \begin{equation}
\mathcal{G}_{\textnormal{mac},i+1/2}^{n,s}
  =
    \frac{1}{2}
    \left[
        \bfG
        \left(
            \bfU_{i+1/2}^{-,n,s},0
        \right)
        +
        \bfG
        \left(
            \bfU_{i+1/2}^{+,n,s},0
        \right)
    \right]
    -
    \frac{a_{i+1/2}^{n,s}}{2}
    \left(
        \bfU_{i+1/2}^{+,n,s}
        -
        \bfU_{i+1/2}^{-,n,s}
    \right)
,
    \label{eq:macroscopic_numerical_flux}
    \end{equation}
where $a_{i+1/2}^{n,s} = \max \left\{ \left|u_{i+1/2}^{-,n,s}\right| + \sqrt{3\theta_{i+1/2}^{-,n,s}} \;,\; \left|u_{i+1/2}^{+,n,s}\right| + \sqrt{3\theta_{i+1/2}^{+,n,s}} \right\}$ is the local Rusanov wave speed used by the macroscopic flux.
The remaining reconstruction, finite-volume, and IMEX details are given in Appendix~\ref{app:kinetic_solver}.

The two numerical fluxes differ in the order in which the kinetic information is reduced and discretized.
The kinetic solver first reconstructs and upwinds the distribution function separately at each discrete velocity, as in \eqref{eq:kinetic_distribution_flux}, and then takes the conserved velocity moments in \eqref{eq:kinetic_numerical_moment_flux}.
The macroscopic solver instead takes the moment system \eqref{eq:kinetic_macro_system} as its starting point and constructs the Rusanov flux \eqref{eq:macroscopic_numerical_flux} directly from the reconstructed macroscopic states.
Consequently, these two discrete approximations need not agree, even when the kinetic heat flux is known exactly.

For an $S$-stage IMEX Runge--Kutta method, define the corresponding stage-weighted transport fluxes as
    \begin{equation}
\mathcal{G}_{\textnormal{kin},i+1/2}^n
  =
    \sum_{s=1}^{S}
      b_s^E
      \mathcal{G}_{\textnormal{kin},i+1/2}^{n,s}
,\qquad
\mathcal{G}_{\textnormal{mac},i+1/2}^n
  =
    \sum_{s=1}^{S}
      b_s^E
      \mathcal{G}_{\textnormal{mac},i+1/2}^{n,s}
,
    \label{eq:stage_weighted_fluxes}
    \end{equation}
where $b_s^E$ are the explicit transport weights.
We define the discrepancy between these fluxes as
    \begin{equation}
\bfR_{\textnormal{eff},i+1/2}^n
  =
    \mathcal{G}_{\textnormal{kin},i+1/2}^n
    -
    \mathcal{G}_{\textnormal{mac},i+1/2}^n
.
    \label{eq:R_eff_definition}
    \end{equation}
Here, $\bfR_{\textnormal{eff}}$ represents the error induced by the discrepancy between the stage-weighted macroscopic and kinetic numerical transport fluxes in \eqref{eq:stage_weighted_fluxes}.
Equivalently, $\bfR_{\textnormal{eff}}$ combines the contribution that would ordinarily enter through the unresolved heat flux with the remaining discretization mismatch between the kinetic and macroscopic numerical transport fluxes.

Since the kinetic data have a negligible macroscopic contribution from the BGK collision term, learning $\bfR_{\textnormal{eff}}$ defined in \eqref{eq:R_eff_definition} is sufficient to capture the discrepancy between the fluxes in the kinetic and macroscopic solvers.
On the other hand, if the kinetic data have a non-negligible macroscopic contribution from the BGK collision term, this contribution can be accounted for through a separate collision-defect correction term in the macroscopic update.
In practice, once the closure term $\bfR_{\textnormal{eff}}$ is learned, the nonlinear discrete-equilibrium solve described in Appendix~\ref{app:kinetic_solver} is no longer needed when solving the macroscopic equation for prediction.

This motivates us to identify this error term as a function of the macroscopic field, $\bfR_{\textnormal{eff}} = \widetilde{\mathcal{R}}(\bfU)$.
Given such a closure operator, $\widetilde{\mathcal{R}}$, and a macroscopic state $\widehat{\bfU}^n$ at time $t_n$, during each prediction step, one corrects the numerical flux as follows:
    \begin{equation}
\widehat{\mathcal{G}}_{i+1/2}^n
  =
    \mathcal{G}_{\textnormal{mac},i+1/2}^n
    \big(
        \widehat{\bfU}^n
    \big)
    +
    \widehat{\bfR}_{\textnormal{eff},i+1/2}^n
,\qquad
\widehat{\bfR}_{\textnormal{eff}}^n
  =
    \widetilde{\mathcal{R}}
    \big(
        \widehat{\bfU}^n
    \big)
,
    \label{eq:learned_effective_flux}
    \end{equation}
Subsequently, the conservative macroscopic update is given as,
    \begin{equation}
\widehat{\bfU}_i^{n+1}
  =
    \widehat{\bfU}_i^n
    -
    \frac{\Delta t}{\Delta x}
    \left(
        \widehat{\mathcal{G}}_{i+1/2}^n
        -
        \widehat{\mathcal{G}}_{i-1/2}^n
    \right)
.
    \label{eq:kinetic_hybrid_rollout}
    \end{equation}

With the closure framework in \eqref{eq:learned_effective_flux}, the remaining modeling question is how to represent the map $\widetilde{\mathcal{R}}$, which we discuss in the following subsection.

\subsection{Closure model representations}
\label{sec:kinetic_representations}

We next consider several representations of the map $\widetilde{\mathcal{R}}$ for the flux correction in \eqref{eq:learned_effective_flux}.
All of the models considered use KRR to approximate this closure map, but they differ in how the spatial state and correction are represented.
We only report results using the Gaussian kernel since we do not find any significant improvement using the DM kernel in this application.

\paragraph{Global Model.}

We first treat the complete macroscopic field and effective flux correction as a single global variable, $\bfU^n,\bfR_{\textnormal{eff}}^n\in\bbR^{3N_x}$.
Because these full-field vectors are high dimensional, PCA is applied separately to the state and correction fields:
Let $\boldsymbol{\xi}_{\bfU}^n\coloneqq\mathcal{P}_{\bfU}(\bfU^n)$ and $\boldsymbol{\xi}_{\bfR}^n\coloneqq \mathcal{P}_{\bfR}(\bfR_{\textnormal{eff}}^n)$ represent the corresponding PCA coordinates, where $\mathcal{P}_{\bfU}$ and $\mathcal{P}_{\bfR}$ denote the PCA projection operators.
Collect the training sets as $\Xi_{\bfU} = \{\boldsymbol{\xi}_{\bfU}^n\}$ and $\Xi_{\bfR} = \{\boldsymbol{\xi}_{\bfR}^n\}$.
The global closure model therefore has the form
    \begin{equation}
\widetilde{\mathcal{R}}({\bfU}^n)
  =
    \mathcal{P}_{\bfR}^{-1}(
      \underbrace{
          A_{\varepsilon,\lambda_{\textnormal{reg}}}^\top
          k_\varepsilon
          (\Xi_{\bfU},\mathcal{P}_{\bfU}(\bfU^n))
      }_{
          \text{KRR on PC components}
      }
    )
,\qquad
A_{\varepsilon,\lambda_{\textnormal{reg}}}
  =
    \left(
        K_\varepsilon
        +
        \lambda_{\textnormal{reg}} I
    \right)^{-1}
    \Xi_{\bfR}
,
    \end{equation}
where $K_\varepsilon$ denotes the Gram matrix induced by data set $\Xi_{\bfU}$ in the PC coordinates.

Preliminary global-model experiments demonstrate accurate predictions for unimodal initial conditions, but also reveal substantial sensitivity to spatial translations and asymmetric bimodal initial data.
Explicitly accounting for periodic translation symmetry can largely remove error due to translations within the unimodal family, but the issue of asymmetric bimodal initial data remains a difficult generalization for the global model.
The local formulation is introduced next to address these difficulties.

\paragraph{Zeroth-order local model.}

While the global approach above is effective for unimodal initial conditions, it performs worse on more challenging cases with asymmetric, bimodal initial conditions.
To overcome this issue, we use spatial localization rather than PCA.
For the present application, consider spatial localization operator, $S_{L,i}$,
    \[
S_{L,i}(\bfU^n)= \bfU_{\textnormal{loc},i}^n
  =
    \left(
        \bfU_{i-L/2+1}^n,\ldots,\bfU_{i+L/2}^n
    \right)
    \in\bbR^{3L}
,
    \]
and set the localization length to be $L=10$.
Subsequently, we learn a KRR model to detect the correction for each interface cell $i+1/2$:
    \[
\bfU_{\textnormal{loc},i}^n
  \overset
      {\widetilde{\mathcal{R}}_{0}}
      {\longmapsto}
    \bfR_{\textnormal{eff},i+1/2}^n
    \in\bbR^3
.
    \]
We note that the local KRR model is a single model attained from training on a set of label data sampled locally in space and at different times, and evaluated at each spatial interface.

\paragraph{First-order local model.}

We also consider a higher-order local closure in which a second KRR model is introduced to correct the error induced by modeling with the zeroth-order model $\widetilde{\mathcal{R}}_0$.
Define
    \[
\bfr_{i+1/2}^n
  =
    \bfR_{\textnormal{eff},i+1/2}^n
    -
    \widetilde{\mathcal{R}}_0
    \left(
        \bfU_{\textnormal{loc},i}^n
    \right)
.
    \]
We then train a KRR model to identify the following map,
    \[
\mathcal{D}\bfU_{\textnormal{loc},i}^n
  \overset{\widetilde{\mathcal{R}}_{1}}{\longmapsto}
    \bfr_{i+1/2}^n
    \in\bbR^3
,
    \]
where
    \[
\mathcal{D}\bfU_{\textnormal{loc},i}^n
  \coloneqq
    \left(
        \frac{\bfU_{i-L/2+2}^n-\bfU_{i-L/2+1}^n}{\Delta x}
        ,\ldots,
        \frac{\bfU_{i+L/2}^n - \bfU_{i+L/2-1}^n}{\Delta x}
    \right)
    \in\bbR^{3(L-1)}
,
    \]
denote the $L-1$ adjacent differences of each of the three conserved variables in the $L$-cell stencil.
With this correction, the first-order local closure takes the form
    \begin{equation}
\widehat{\bfR}_{\textnormal{eff},i+1/2}^n
  =
    \widetilde{\mathcal{R}}_{0}
    \left(
        \widehat{\bfU}_{\textnormal{loc},i}^n
    \right)
    +
    \widetilde{\mathcal{R}}_{1}
    \left(
        \mathcal{D}
        \widehat{\bfU}_{\textnormal{loc},i}^n
    \right)
.
    \label{eq:kinetic_derivative_residual}
    \end{equation}
At every macroscopic time step, both terms in \eqref{eq:kinetic_derivative_residual} are reevaluated from the current recursively predicted state $\widehat{\bfU}^n$; the residual variable is not advanced independently from $t^n$ to $t^{n+1}$.

\

The global, zeroth-order local, and first-order local representations therefore use the same effective closure target but impose substantially different geometries on the KRR regression problem.
The training sets and rollout-based evaluation procedure used to compare these representations are described next.

\subsection{Training data and evaluation metric}
\label{sec:kinetic_evaluation}

All kinetic trajectories are initialized from the family in \eqref{eq:macroscopic_ic_fields}.
We organize these trajectories into two data types according to the structure of the initial conditions.
Data Type I consists of unimodal initial conditions, with $\mu=\eta\Delta x$.
Data type II consists of bimodal initial conditions with modes at $\mu_{1,2} = \phi \pm \frac{d}{2}$.
The parameters varied within each data type are summarized in Table~\ref{table:kinetic_data_types}, and specific parameter values and their assignments to training, validation, and testing datasets are given in Appendix~\ref{app:kinetic_data}.

\begin{table}[htbp]
    \centering
    \small
    \begin{tabular}{cll}
        \toprule
        Data type
        & Phase shift
        & Varying parameters
        \\
        \midrule

        I
        & $\mu=\eta\Delta x$
        & $A,\ \sigma,\ \eta$
        \\

        II
        & $\mu_{1,2} = \phi \pm \frac{d}{2}$
        & $A_1,A_2,\sigma_1,\sigma_2,$ $d$, $\phi$
        \\

        \bottomrule
    \end{tabular}
    \caption{
        The two physical initial condition data types used in the kinetic experiments.
    }
    \label{table:kinetic_data_types}
\end{table}

We will specify the types of training data being used in constructing the Gram matrix of KRR in Section~\ref{sec:kinetic_results}.
We now specify the evaluation metrics to be used for validation and testing.

Denote $\mathcal{S}_{\text{I}}$ and $\mathcal{S}_{\text{II}}$ as collections of trajectories with unimodal (Type I) and bimodal (Type II) initial conditions, respectively.
Given the $\tau^{\textnormal{th}}$ member of $p_{\text{true},\tau} \in S_j$, where $p \in \{\rho,u,\theta,E\}$ is a spatio-temporal field of temporal length $N_\tau$ and $j\in\{\text{I},\text{II}\}$, we define the prediction at time $t^n$ from the closure model given initial condition $p_{\text{true},\tau}(0)$ as $p_{\text{pred},\tau}(n)$, which is obtained by iterating the closure model \eqref{eq:learned_effective_flux}--\eqref{eq:kinetic_hybrid_rollout}.

We evaluate the prediction accuracy with the following metric,
    \begin{equation}
\mathcal{E}(\mathcal{S}_j)
  \coloneqq
    \frac{1}{4|\mathcal{S}_j|}
    \sum_{\tau \in \mathcal{S}_j}
      \sum_{p\in \{\rho,u,\theta,E\}}
        e_{\tau,p}(N_\tau)
,
    \end{equation}
where
    \[
e_{\tau,p}(n)
  =
    \frac{
        \|
          p_{\textnormal{true},\tau}(n)
          -
          p_{\textnormal{pred},\tau}(n)
        \|_{2,h}
    }{
        \|
          p_{\textnormal{true},\tau}(n)
        \|_{2,h}
    }
,\qquad
\|p\|_{2,h}
  =
    \left(
        \Delta x\sum_i |p_i|^2
    \right)^{1/2}
,\qquad
    p\in\{\rho,u,\theta,E\}
,
    \]
is a normalized spatial $L^2-$norm.

The above metric is used during the validation procedure to select model hyperparameters $\varepsilon$ and $\lambda_{\textnormal{reg}}$ in Section~\ref{sec:DMKRR}.
In particular, let $\mathcal{V}_{\text{I}}$ and $\mathcal{V}_{\text{II}}$ denote the sets of validation trajectories constructed from unimodal and bimodal initial conditions respectively.
We validate the KRR hyperparameters by minimizing, $\mathcal{E}(\mathcal{V}_{\text{I}})+ \mathcal{E}(\mathcal{V}_{\text{II}})$.
Once the model is trained, the same scoring rule is applied to the test families $\mathcal{T}_{\text{I}}$ and $\mathcal{T}_{\text{II}}$.
We report the testing errors $\mathcal{E}_1 \coloneqq \mathcal{E}(\mathcal{T}_{\text{I}})$ and $\mathcal{E}_2 \coloneqq \mathcal{E}(\mathcal{T}_{\text{II}})$ in Table~\ref{table:kinetic_model_comparison}.

\subsection{Results}
\label{sec:kinetic_results}

The numerical results are summarized in Table~\ref{table:kinetic_model_comparison}.
The  KRR centers column in the table corresponds to the size $N$ of training datasets used to fit the regression models.
For the first-order local closure model, we report $(N_0,N_1)$, where $N_0$ and $N_1$ are the center counts used in training $\widetilde{\mathcal{R}}_0$ and $\widetilde{\mathcal{R}}_1$, respectively.

The global models trained using different combinations of training data types attain small errors ($\mathcal{E}_1$) in predicting test data with unimodal initial conditions.
The errors are substantially larger for bimodal initial conditions, although adding Type II trajectories to the global training set considerably reduces $\mathcal{E}_2$.

Spatial localization produces a further substantial reduction in the bimodal error, $\mathcal{E}_2$.
The zeroth-order local and first-order local models exhibit different accuracy tradeoffs across the two test families: the first-order residual model is more accurate on the unimodal testing data, whereas the zeroth-order local model is slightly more accurate on the bimodal family under mixed training.
Averaged across the two test families, the first-order local model gives the smallest error among the three representations.

\begin{table}[htbp]
    \centering
    \small
    \begin{tabular}{lcccc}
        \toprule
        Model
        & Training data
        & KRR centers
        & $\mathcal{E}_1$
        & $\mathcal{E}_2$
        \\
        \midrule

        Global
        & I
        & $8208$
        & $0.001285 \pm 0.000200$
        & $0.028925 \pm 0.005158$
        \\

        Global
        & I+II
        & $8208$
        & $0.001005 \pm 0.000234$
        & $0.008638 \pm 0.007220$
        \\

        \midrule

        zeroth-order local
        & I
        & $4032$
        & $0.001455 \pm 0.000252$
        & $0.006088 \pm 0.001457$
        \\

        zeroth-order local
        & I+II
        & $4032$
        & $0.001891 \pm 0.000341$
        & $0.003178 \pm 0.000675$
        \\

        \midrule

        first-order local
        & $(\textnormal{I},\,\textnormal{I})$
        & $(4032,8064)$
        & $0.000665 \pm 0.000096$
        & $0.008098 \pm 0.002001$
        \\

        first-order local
        & $(\textnormal{I},\,\textnormal{I+II})$
        & $(4032,8064)$
        & $0.000793 \pm 0.000105$
        & $0.003768 \pm 0.001161$
        \\

        \bottomrule
    \end{tabular}
    \caption{
        Comparison of the kinetic closure models under different closure models and training-data compositions.
        The columns $\mathcal{E}_1$ and $\mathcal{E}_2$ report the average prediction errors for the unimodal and bimodal testing data, respectively (lower values are better).
        For the first-order local model, the ordered pair in the training-data column gives the trajectory types used to train $(\widetilde{\mathcal{R}}_0,\widetilde{\mathcal{R}}_1)$, and the number of training centers is reported as $(N_0,N_1)$.
        The zeroth-order model $\widetilde{\mathcal{R}}_0$ is held fixed between the two first-order comparisons.
    }
    \label{table:kinetic_model_comparison}
\end{table}

\begin{figure}[htbp]
    \centering
    \includegraphics[
        width=0.92\textwidth
    ]{
        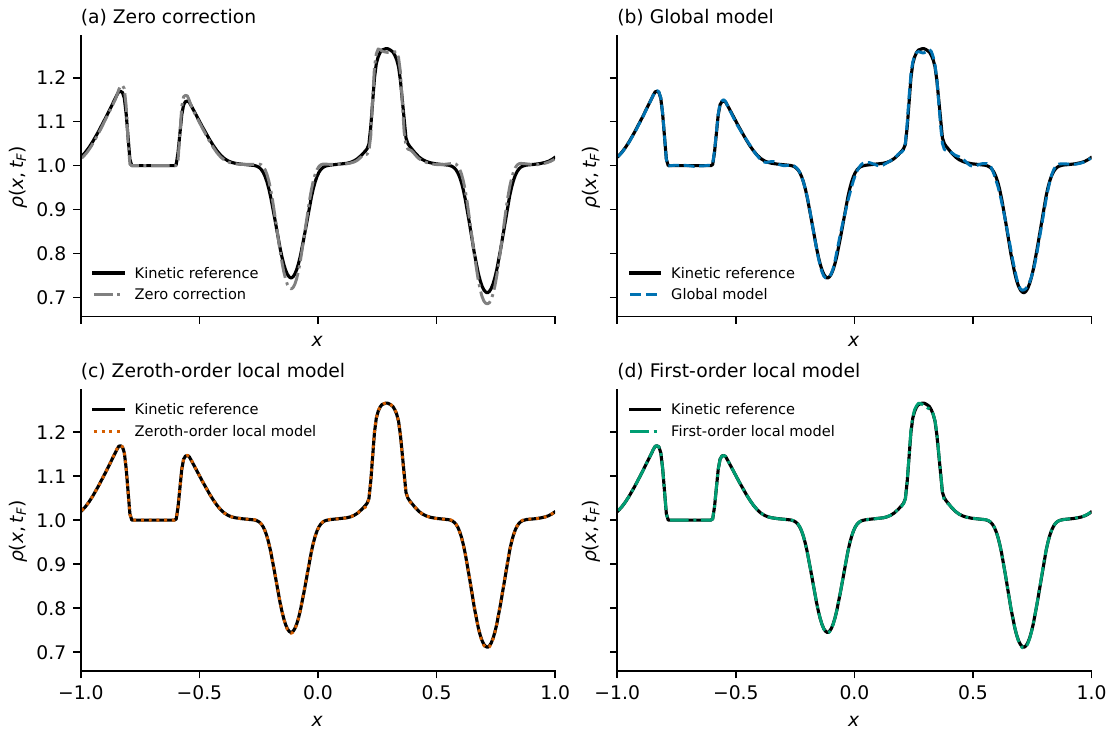
    }
    \caption{
        Final-time density profiles for a Type II validation trajectory.
        Each panel compares the kinetic reference with one macroscopic model: the zero-correction baseline, the global PCA--KRR closure, the zeroth-order local closure, and the first-order local residual closure.
    }
    \label{fig:kinetic_pairwise_rho}
\end{figure}

\begin{figure}[t]
    \centering
    \includegraphics[
        width=\textwidth
    ]{
        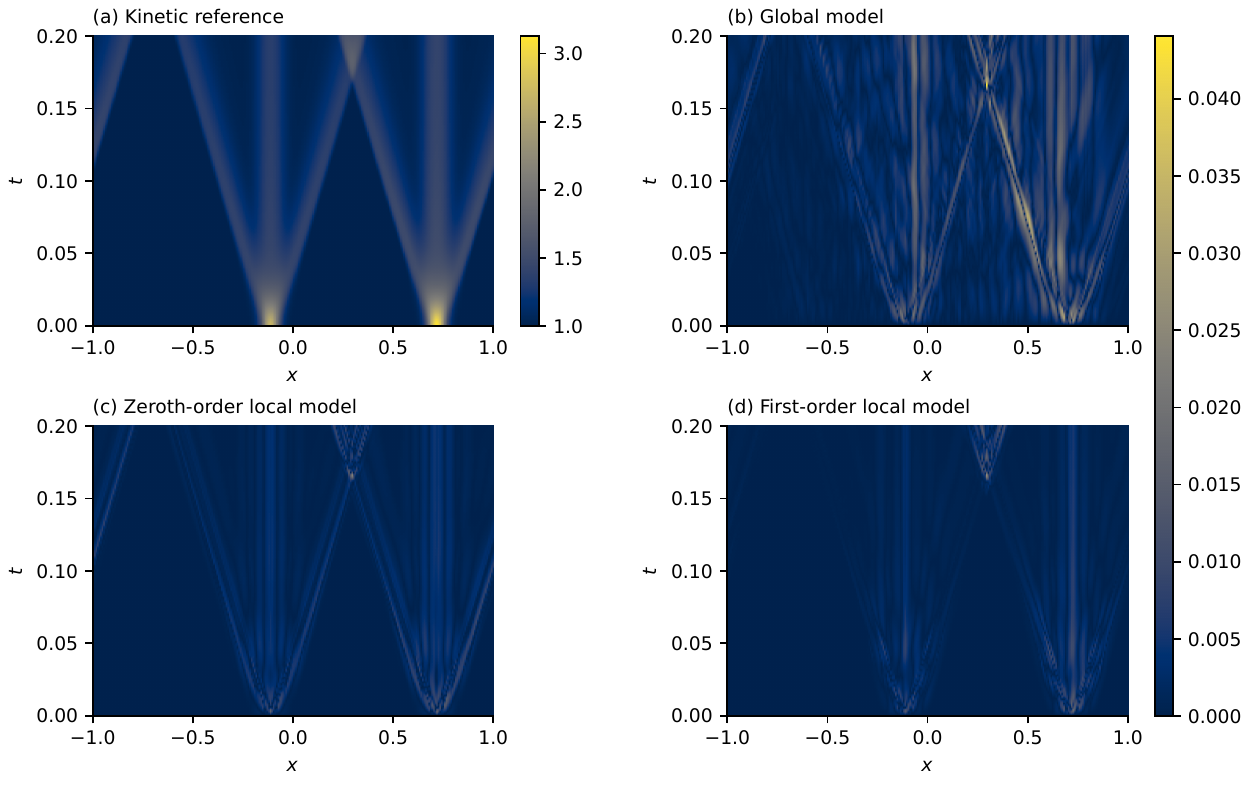
    }
    \caption{
        Space--time evolution of temperature for the same Type II validation trajectory as in Figure~\ref{fig:kinetic_pairwise_rho}.
        Panel (a) shows the kinetic reference temperature $\theta$, while panels (b)--(d) show the absolute temperature errors, for the global model, the zeroth-order local, and the first-order local model, respectively.
        The three error panels use a common linear color scale so that their magnitudes are directly comparable.
        The global model can reproduce the final macroscopic profile reasonably well while accumulating substantially larger structured transient errors during recursive rollout; both local models substantially reduce the structured rollout error and track the kinetic dynamics more faithfully throughout the trajectory.
    }
    \label{fig:kinetic_theta_spacetime}
\end{figure}

The results demonstrated in Table~\ref{table:kinetic_model_comparison} are also visible in the predicted macroscopic fields.
Figure~\ref{fig:kinetic_pairwise_rho} compares the final-time density profile of each model separately with the kinetic reference, demonstrating differences in profile shape and spurious oscillations.
The corresponding final-time velocity, temperature, and energy profiles are shown in Figures~\ref{fig:kinetic_additional_u}--\ref{fig:kinetic_additional_E} of Appendix~\ref{app:kinetic_additional_figures}.
The corresponding space--time temperature errors in Figure~\ref{fig:kinetic_theta_spacetime} show when these discrepancies develop.

\section{Conclusion}
\label{sec:conclusion}

We have developed a closure modeling framework for identifying missing components of dynamical systems using Kernel Ridge Regression (KRR).
We examined the proposed approach on two classes of closure problems.
The first arises in ODE and PDE settings, where the missing components take the form of difference equations.
The second arises from the moment closure problem in kinetic equations, where the closure model is formulated as a system of algebraic equations.

For the first class of problems, we derived an error bound for the proposed framework in an ODE setting and identified the contributions to the total error from the ODE solver, the approximation of the unresolved scales, and the interpolation of the unresolved-scale contributions needed to couple them to the resolved numerical solver.
We numerically tested the framework on two benchmark problems: the Lorenz-63 system and the Kuramoto–Sivashinsky equation.
In both examples, the proposed KRR framework produced accurate long-horizon predictions, with substantial improvements over the LSTM model.
Of the two kernels considered, the diffusion maps (DM) kernel yielded slightly more accurate predictions.

For the second class of problems, we considered moment closure for a one-dimensional kinetic equation.
In this setting, our closure approach models the discrepancy between the numerical transport flux constructed by a computationally expensive kinetic solver that preserves the discrete collision invariants using the conservative discrete Maxwellian of Mieussens~\cite{Mieussens2000}, whose parameters are obtained using a nonlinear solver that matches the discrete conserved moments to the solver tolerance.
Our numerical experiments show that the predictive performance of the closure model is sensitive to the choice of model configuration.
Specifically, we considered three flux-correction models:
(1) a global model using KRR in PCA coordinates;
(2) a zeroth-order spatially local model using local resolved states as inputs;
and (3) a first-order spatially local residual model that augments the zeroth-order model with an additional learned correction using spatial differences of the resolved states.
Numerically, we found that the global model achieves high accuracy for unimodal initial conditions but is less accurate for bimodal initial conditions.
The zeroth-order local model improves the prediction accuracy, while the first-order local model achieves the highest accuracy on average across the two test families.
These results indicate that spatial localization is the principal source of the improved robustness observed in the kinetic experiments.
Practically, the proposed closure scheme allows the coupled system to be solved using only the corrected macroscopic solver, avoiding the expensive kinetic computation and the associated nonlinear discrete-equilibrium scheme \cite{Mieussens2000} required to maintain discrete conservation in the kinetic solver.

Future work will investigate extensions of the framework to noisy data using the weak formulation recently proposed in \cite{kreider2026learning}.
An open question arising from this study, particularly from the empirical results for the Kuramoto–Sivashinsky equation, is why non-Markovian models do not appear to be essential within the proposed closure framework.
The numerical results in the kinetic problem motivate extensions to more challenging kinetic regimes, including nonsmooth flows, higher-dimensional problems, and more general collision models.

\section*{Code availability}

Code, experiment configurations, and instructions for reproducing the KRR experiments are available at
\url{https://github.com/evanhabbershaw/learning-closure-krr}.

\section*{Acknowledgments}

The research of J.H. was partially supported by the NSF grant DMS-2505605.
E.H. was partially supported by the Center for Computational Mathematics and Applications at the Pennsylvania State University.
S.L. was partially supported by the startup fund from Texas Tech University.

\bibliographystyle{plain}
\bibliography{refs}

@article{McKayBeckmanConover1979,
  author  = {McKay, M. D. and Beckman, R. J. and Conover, W. J.},
  title   = {A Comparison of Three Methods for Selecting Values of Input Variables in the Analysis of Output from a Computer Code},
  journal = {Technometrics},
  volume  = {21},
  number  = {2},
  pages   = {239--245},
  year    = {1979},
  doi     = {10.1080/00401706.1979.10489755}
}

@article{MorrisMitchell1995,
  author  = {Morris, Max D. and Mitchell, Toby J.},
  title   = {Exploratory Designs for Computational Experiments},
  journal = {Journal of Statistical Planning and Inference},
  volume  = {43},
  number  = {3},
  pages   = {381--402},
  year    = {1995},
  doi     = {10.1016/0378-3758(94)00035-T}
}

@article{kingma2014adam,
  title   = {Adam: A method for stochastic optimization},
  author  = {Kingma, Diederik P and Ba, Jimmy},
  journal = {arXiv preprint arXiv:1412.6980},
  year    = {2014}
}

@article{hamzi2021learning,
  title     = {Learning dynamical systems from data: a simple cross-validation perspective, part I: parametric kernel flows},
  author    = {Hamzi, Boumediene and Owhadi, Houman},
  journal   = {Physica D: Nonlinear Phenomena},
  volume    = {421},
  pages     = {132817},
  year      = {2021},
  publisher = {Elsevier}
}

@article{darcy2025learning,
  title     = {Learning dynamical systems from data: a simple cross-validation perspective, part {II}: nonparametric kernel flows},
  author    = {Darcy, Matthieu and Hamzi, Boumediene and Susiluoto, Jouni and Braverman, Amy and Owhadi, Houman},
  journal   = {Physica D: Nonlinear Phenomena},
  volume    = {476},
  pages     = {134641},
  year      = {2025},
  publisher = {Elsevier}
}

@article{levine2022framework,
  title   = {A framework for machine learning of model error in dynamical systems},
  author  = {Levine, Matthew and Stuart, Andrew},
  journal = {Communications of the American Mathematical Society},
  volume  = {2},
  number  = {7},
  year    = {2022}
}

@article{crommelin2008subgrid,
  title   = {Subgrid-scale parameterization with conditional {Markov} chains},
  author  = {Crommelin, Daan and Vanden-Eijnden, Eric},
  journal = {Journal of the Atmospheric Sciences},
  volume  = {65},
  number  = {8},
  pages   = {2661--2675},
  year    = {2008}
}

@article{lu2016comparison,
  title     = {Comparison of continuous and discrete-time data-based modeling for hypoelliptic systems},
  author    = {Lu, Fei and Lin, Kevin and Chorin, Alexandre},
  journal   = {Communications in Applied Mathematics and Computational Science},
  volume    = {11},
  number    = {2},
  pages     = {187--216},
  year      = {2016},
  publisher = {Mathematical Sciences Publishers}
}

@article{mh:13,
  Author  = {Majda, A.J. and Harlim, J.},
  Journal = {Nonlinearity},
  Pages   = {201-217},
  Title   = {{Physics constrained nonlinear regression models for time series.}},
  Volume  = {26},
  Year    = {2013}
}

@article{hmm:14,
  Author  = {Harlim, J. and Mahdi, A. and Majda, A.J.},
  Journal = {J. Comput. Phys.},
  Pages   = {782-812},
  Title   = {{An ensemble Kalman filter for statistical estimation of physics constrained nonlinear regression models}},
  Volume  = {257, Part A},
  Year    = {2014}
}

@article{bh:14,
  Author  = {Berry, T. and Harlim, J.},
  Journal = {Proc. Roy. Soc. A 20140168},
  Title   = {{Linear Theory for Filtering Nonlinear Multiscale Systems with Model Error}},
  Year    = {2014}
}

@article{lu2017accounting,
  title   = {Accounting for Model Error from Unresolved Scales in Ensemble Kalman Filters by Stochastic Parameterization},
  author  = {Lu, Fei and Tu, Xuemin and Chorin, Alexandre J},
  journal = {Monthly Weather Review},
  volume  = {145},
  number  = {9},
  pages   = {3709--3723},
  year    = {2017}
}

@article{kwasniok2012data,
  title     = {Data-based stochastic subgrid-scale parametrization: an approach using cluster-weighted modelling},
  author    = {Kwasniok, Frank},
  journal   = {Philosophical Transactions of the Royal Society A: Mathematical, Physical and Engineering Sciences},
  volume    = {370},
  number    = {1962},
  pages     = {1061--1086},
  year      = {2012},
  publisher = {The Royal Society Publishing}
}

@article{Mieussens2000,
  author  = {Mieussens, Luc},
  title   = {Discrete Velocity Model and Implicit Scheme for the
             {BGK} Equation of Rarefied Gas Dynamics},
  journal = {Mathematical Models and Methods in Applied Sciences},
  volume  = {10},
  number  = {8},
  pages   = {1121--1149},
  year    = {2000}
}

@article{HanMaMaE2019,
  author  = {Han, Jiequn and Ma, Chao and Ma, Zheng and E, Weinan},
  title   = {Uniformly Accurate Machine Learning-Based Hydrodynamic Models
             for Kinetic Equations},
  journal = {Proceedings of the National Academy of Sciences},
  volume  = {116},
  number  = {44},
  pages   = {21983--21991},
  year    = {2019},
  doi     = {10.1073/pnas.1909854116}
}

@article{LiDongWang2023,
  author  = {Li, Zhengyi and Dong, Bin and Wang, Yanli},
  title   = {Learning Invariance Preserving Moment Closure Model for
             Boltzmann--BGK Equation},
  journal = {Communications in Mathematics and Statistics},
  volume  = {11},
  number  = {1},
  pages   = {59--101},
  year    = {2023},
  doi     = {10.1007/s40304-022-00331-5}
}

@article{HuangDongWang2025,
  author  = {Huang, Ziyu and Dong, Chuanfei and Wang, Liang},
  title   = {Machine-Learning Heat Flux Closure for Multi-Moment Fluid
             Modeling of Nonlinear Landau Damping},
  journal = {Proceedings of the National Academy of Sciences},
  volume  = {122},
  number  = {11},
  pages   = {e2419073122},
  year    = {2025},
  doi     = {10.1073/pnas.2419073122}
}

@article{IngelstenEtAl2025,
  author  = {Ingelsten, Emil R. and McGrae-Menge, Madox C.
             and Alves, E. Paulo and Pusztai, Istvan},
  title   = {Data-Driven Discovery of a Heat Flux Closure for
             Electrostatic Plasma Phenomena},
  journal = {Journal of Plasma Physics},
  volume  = {91},
  number  = {2},
  pages   = {E64},
  year    = {2025},
  doi     = {10.1017/S0022377825000285}
}

@article{BurlesCamporeale2026,
  author  = {Burles, Samuel and Camporeale, Enrico},
  title   = {The Machine Learning Approach to Moment Closure Relations
             for Plasma: A Review},
  journal = {Journal of Plasma Physics},
  volume  = {92},
  number  = {4},
  pages   = {E100},
  year    = {2026},
  doi     = {10.1017/S0022377826102013}
}

@article{BhatnagarGrossKrook1954,
  author  = {Bhatnagar, P. L. and Gross, E. P. and Krook, M.},
  title   = {A Model for Collision Processes in Gases. {I}. Small Amplitude
             Processes in Charged and Neutral One-Component Systems},
  journal = {Physical Review},
  volume  = {94},
  number  = {3},
  pages   = {511--525},
  year    = {1954},
  doi     = {10.1103/PhysRev.94.511}
}

@book{CercignaniIllnerPulvirenti1994,
  author    = {Cercignani, Carlo and Illner, Reinhard and Pulvirenti, Mario},
  title     = {The Mathematical Theory of Dilute Gases},
  series    = {Applied Mathematical Sciences},
  volume    = {106},
  publisher = {Springer},
  year      = {1994},
  doi       = {10.1007/978-1-4419-8524-8}
}

@article{Levermore1996,
  author  = {Levermore, C. David},
  title   = {Moment Closure Hierarchies for Kinetic Theories},
  journal = {Journal of Statistical Physics},
  volume  = {83},
  pages   = {1021--1065},
  year    = {1996},
  doi     = {10.1007/BF02179552}
}

@article{PareschiRusso2005,
  author  = {Pareschi, Lorenzo and Russo, Giovanni},
  title   = {Implicit-Explicit Runge-Kutta Schemes and Applications to
             Hyperbolic Systems with Relaxation},
  journal = {Journal of Scientific Computing},
  volume  = {25},
  number  = {1},
  pages   = {129--155},
  year    = {2005},
  doi     = {10.1007/s10915-004-4636-4}
}

@article{HuangChengChristliebRoberts2022,
  author  = {Huang, Juntao and Cheng, Yingda and Christlieb, Andrew J.
             and Roberts, Luke F.},
  title   = {Machine Learning Moment Closure Models for the Radiative
             Transfer Equation {I}: Directly Learning a Gradient Based Closure},
  journal = {Journal of Computational Physics},
  volume  = {453},
  pages   = {110941},
  year    = {2022},
  doi     = {10.1016/j.jcp.2022.110941}
}

@article{PorteousLaiuHauck2023,
  author  = {Porteous, William A. and Laiu, M. Paul and Hauck, Cory D.},
  title   = {Data-Driven, Structure-Preserving Approximations to Entropy-Based
             Moment Closures for Kinetic Equations},
  journal = {Communications in Mathematical Sciences},
  volume  = {21},
  number  = {4},
  pages   = {885--913},
  year    = {2023},
  doi     = {10.4310/CMS.2023.v21.n4.a1}
}

@article{ChristliebDingHuangKrupansky2025,
  author  = {Christlieb, Andrew J. and Ding, Mingchang and Huang, Juntao
             and Krupansky, Nicholas A.},
  title   = {Hyperbolic Machine Learning Moment Closures for the {BGK} Equations},
  journal = {Multiscale Modeling \& Simulation},
  volume  = {23},
  number  = {1},
  pages   = {187--217},
  year    = {2025},
  doi     = {10.1137/24M1629377}
}

@book{Toro2009,
    author    = {Eleuterio F. Toro},
    title     = {Riemann Solvers and Numerical Methods for Fluid Dynamics:
                 A Practical Introduction},
    edition   = {3},
    publisher = {Springer},
    address   = {Berlin, Heidelberg},
    year      = {2009},
    doi       = {10.1007/b79761}
}

@article{Rusanov1962,
    author  = {V. V. Rusanov},
    title   = {The Calculation of the Interaction of Non-Stationary
               Shock Waves and Obstacles},
    journal = {USSR Computational Mathematics and Mathematical Physics},
    volume  = {1},
    number  = {2},
    pages   = {304--320},
    year    = {1962},
    doi     = {10.1016/0041-5553(62)90062-9}
}

@article{harlim2026diffusion,
  title   = {Diffusion Maps Kernel Ridge Regression},
  author  = {Harlim, John and Huang, Daning and Song, Jiwoo},
  journal = {arXiv preprint arXiv:2608.09623},
  year    = {2026},
  url     = {https://arxiv.org/abs/2608.09623}
}

@article{kreider2026learning,
  title   = {Learning dynamical systems from noisy data with Weak-form Kernel Ridge Regression},
  author  = {Kreider, Max and Harlim, John and Huang, Daning},
  journal = {arXiv preprint arXiv:2607.00257},
  year    = {2026}
}

@article{mandal2025learning,
  title     = {Learning dynamical systems with hit-and-run random feature maps},
  author    = {Mandal, Pinak and Gottwald, Georg A},
  journal   = {Nature Communications},
  volume    = {16},
  number    = {1},
  pages     = {5961},
  year      = {2025},
  publisher = {Nature Publishing Group UK London}
}

@article{coifman2008graph,
  title     = {Graph {Laplacian} tomography from unknown random projections},
  author    = {Coifman, Ronald R and Shkolnisky, Yoel and Sigworth, Fred J and Singer, Amit},
  journal   = {IEEE Transactions on Image Processing},
  volume    = {17},
  number    = {10},
  pages     = {1891--1899},
  year      = {2008},
  publisher = {IEEE}
}

@article{harlim_2021,
  title   = {Machine learning for prediction with missing dynamics},
  journal = {Journal of Computational Physics},
  volume  = {428},
  pages   = {109922},
  year    = {2021},
  issn    = {0021-9991},
  doi     = {10.1016/j.jcp.2020.109922},
  url     = {https://www.sciencedirect.com/science/article/pii/S0021999120306963},
  author  = {John Harlim and Shixiao W. Jiang and Senwei Liang and Haizhao Yang}
}

@article{song2025_KRR_DM,
  title   = {Learning solution operator of dynamical systems with diffusion maps kernel ridge regression},
  author  = {Jiwoo Song and Daning Huang and John Harlim},
  journal = {arXiv preprint 2512.17203},
  year    = {2025},
  url     = {https://arxiv.org/abs/2512.17203}
}

@techreport{habbershaw2025year4,
  author      = {Habbershaw, Evan and Wise, Steven M.},
  title       = {{Year-4 Progress Report on Numerical Methods for BGK-Type Kinetic Equations}},
  institution = {University of Tennessee, Knoxville},
  type        = {Mathematics Report},
  number      = {14},
  year        = {2025},
  url         = {https://trace.tennessee.edu/handle/20.500.14382/48810}
}

\appendix

\section{Rollout-aware model selection}
\label{app:model_selection}

In this appendix, we briefly review the validation strategy reported in \cite{song2025_KRR_DM} for selecting the KRR hyperparameters: the kernel bandwidth $\varepsilon$ and the ridge-regularization parameter $\lambda_{\mathrm{reg}}$.
We begin by obtaining initial guesses for $(\varepsilon,\lambda_{\mathrm{reg}})$ using the procedure in \cite{song2025_KRR_DM}, which is motivated by the intrinsic-dimension estimation algorithm proposed in \cite{coifman2008graph}.

Define the Gaussian kernel,
    \[
\rho(\bfz,\bfz';\eta)
  =
    \exp
    \left(
        -\|\bfz-\bfz'\|^2/(L^2\eta)
    \right)
,
    \]
where $L$ is the maximum of pair-wise L2 distances between elements of the dataset and $\eta$ is the bandwidth parameter.

Given $\{\bfz_i\}_{i=1}^N$ assumed to lie on an unknown manifold of intrinsic dimension $d$, we can estimate $d$ as follows.
Define the sum
    \[
S(\eta)
  =
    \frac{1}{N^2}\sum_{i,j=1}^N \rho(\bfz_i,\bfz_j;\eta)
.
    \]
One can show \cite{coifman2008graph} that $S(\eta)\propto \eta^{d/2}$, which leads to the following estimation of the intrinsic dimension.
Define
    \[
\tilde{d}(\eta)
  =
    \frac{2\mathrm{d}\log(S(\eta))}{\mathrm{d}\log(\eta)}
,
    \]
which can be estimated using a finite difference over a sweep of $\eta$-values.
We empirically found that
    \[
d^*
  =
    \tilde{d}(\eta^*)
\qquad\text{where}\qquad
\eta^*
  =
    \arg\max_\eta \tilde{d}(\eta)
,
    \]
is a good estimate for the intrinsic dimension $d$.

From $\eta^*$, we empirically choose
    \[
\varepsilon^*
  =
    \frac{5}{2d_z}(L^2\eta^*)^{1/d^*}
.
    \]
For chaotic problems, in which the forward invariant set is not smooth, we still use the above algorithm and choose
    \[
\varepsilon^*
  =
    250 L^2\eta^*
.
    \]
Then, using $\varepsilon^*$, construct the symmetric RBF kernel matrix.
Given the minimum eigenvalue of the kernel matrix is $\lambda_{\min}$, choose the reference regularization factor $\lambda_\text{reg}^*=|\lambda_{\min}|$.
The range is chosen to be $H=[\varepsilon^*\Delta \epsilon, \varepsilon^*/\Delta \varepsilon]\times [\lambda_\text{reg}^*, \lambda_\text{reg}^*/\Delta \lambda_\text{reg}]$, where, unless otherwise stated, $\Delta \varepsilon=10^{-2}, \Delta \lambda_\text{reg}=10^{-4}$.
Validation is then employed either by taking randomly sampled parameters or equally log-spaced parameters from $H$.
For each parameter, we compute the VPT as defined in \eqref{eq:L63VPT} on validation trajectories of length $N_V$.
Finally, we select the parameter pair $(\varepsilon,\lambda_{\text{reg}})$ that gives the best VPT.

\section{Proof of Theorem~\ref{thm1}}
\label{app:proof_of_theorem}

Let
    \[
\widehat{\bfx}_{n+1} - \bfx(t_{n+1})
  =
    \big(
      \widehat{\bfx}_{n+1} - \widetilde{\bfx}(t_{n+1})
    \big)
    +
    \big(
      \widetilde{\bfx}(t_{n+1}) - \bfx(t_{n+1})
    \big)
,
    \]
where
    \[
\widetilde{\bfx}'
  =
    f(\widetilde{\bfx},\widehat{y}(t))
,\qquad
t
  \in
    [t_n,t_{n+1}]
,\qquad
\widetilde{\bfx}(t_n)
  =
    \widehat{\bfx}_n
.
    \]
Taking norms gives
    \begin{equation}
e_{n+1}^x
  \coloneqq
    \| \widehat{\bfx}_{n+1} -\bfx(t_{n+1}) \|
  \leq
    \| \widehat{\bfx}_{n+1} - \widetilde{\bfx}(t_{n+1}) \|
    +
    \|\widetilde{\bfx}(t_{n+1}) - \bfx(t_{n+1})\|
.
    \label{eq:pf1}
    \end{equation}
The first term in \eqref{eq:pf1} satisfies
    \[
\|\widehat{\bfx}_{n+1} - \widetilde{\bfx}(t_{n+1})\|
  \leq
    C_0 (\Delta t)^{q+1}
,
    \]
when order-$q$ ODE solver is used.

For any $t \in [t_n,t_{n+1}]$
    \begin{align}
\|\widehat{\bfy}(t) - \bfy(t) \|
  &=
    \|
      \mathcal{I}(\widehat{\bfy}_{n-p+1},\ldots,\widehat{\bfy}_{n+1})(t)
      -
      \bfy(t)
    \|
    \notag
    \\
  &\leq
    \|
      \mathcal{I}(\widehat{\bfy}_{n-p+1},\ldots,\widehat{\bfy}_{n+1})(t)
      -
      \mathcal{I}(\bfy_{n-p+1},\ldots,\bfy_{n+1})(t)
    \|
    \notag
    \\
  &\quad+
    \|
      \mathcal{I}(\bfy_{n-p+1},\ldots,\bfy_{n+1})(t) - \bfy(t)
    \|
    \notag
    \\
  &\leq
    \|\mathcal{I}\|
    \max_{j=0}^p
      \|
        \widehat{\bfy}_{n-j+1} - \bfy_{n-j+1}
      \|
      +
      C_3 (\Delta t)^{p+1}
    \notag
    \\
  &\leq
    C_y
    \left(
      \delta(\epsilon,N) +  (\Delta t)^{p+1}
    \right)
,
    \label{eq:pf2}
    \end{align}
where the two terms in the first inequality above correspond to errors between interpolating the approximate data $\{\widehat{\bfy}_{n-j+1}\}_{j=0,\ldots,p}$ and the true data $\{\bfy_{n-j+1}\}_{j=0,\ldots,p}$ and the interpolation error.
In the second inequality we used the assumption that $\mathcal{I}$ is a stable interpolator with the interpolation error rate $p+1$.
The last inequality corresponds to the error induced by the approximate KRR model.

Define
    \[
\bfw(t)
  \coloneqq
    \widetilde{\bfx}(t) - \bfx(t)
,\qquad
t
  \in
    [t_n,t_{n+1}]
.
    \]
Since $\widetilde{\bfx}(t_n) = \widehat{\bfx}_n$, we have $\|\bfw(t_n)\| = e_n^x$.
One can write
    \[
\bfw_n'(t)
  =
    \widetilde{\bfx}'(t) - \bfx'(t)
  =
    f(\widetilde{\bfx}(t),\widehat{\bfy}(t))
    -
    f(\bfx(t),\bfy(t))
.
    \]
Integrating from $t_n$ to $t\in[t_n,t_{n+1}]$, we obtain
    \[
\int_{t_n}^t \bfw'(s)\,\diff s
  =
    \bfw(t) - \bfw(t_n)
  =
    \int_{t_n}^t
    [
      f(\widetilde{\bfx}_n(s),\widehat{\bfy}(s)) - f(\bfx(s),\bfy(s))
    ]
    \,\diff s
    \]
Rearranging, taking norms, and using triangle inequality and the Lipschitz continuity of $f$,
    \begin{align}
\|\bfw(t)\|
  &\leq
    \|\bfw(t_n)\|
    +
    \int_{t_n}^t
    \Big(
      L_x\|\bfw(s)\|
      +
      L_y\|\widehat{\bfy}(s) - \bfy(s)\|
    \Big)
    \,\diff s
    \notag
    \\
  &\leq
    e_n^x
    +
    L_x
    \int_{t_n}^t
      \|\bfw(s)\|
      \,\diff s
    +
    C_y L_y \Delta t
    \left(
        \delta(\epsilon,N) + (\Delta t)^{p+1}
    \right)
,
    \end{align}
where we have used \eqref{eq:pf2}.
By Gr\"onwall inequality and standard inequality bounds, we have
    \[
\|\bfw(t_{n+1})\|
  \leq
    (1+C_1 \Delta t)e_n^x
    +
    C_2
    \left(
        \Delta t \delta(\epsilon,N)
        +
        (\Delta t)^{p+1}
    \right)
,
    \]
where $C_2 = C_yL_y e^{L_x T}$.
Inserting this into \eqref{eq:pf1} and taking $C = \max\{C_0,C_1,C_2\}$, we have
\[
e_{n+1}^x
  \leq
    (1+C \Delta t)e_n^x
    +
    C\Delta t
    \left(
        \delta(\epsilon,N)
        +
        (\Delta t)^{p+1}
        +
        (\Delta t)^q
    \right)
.
    \]
Using a standard geometric series argument and taking the maximum over $0\leq n \leq N_T$ yields
    \[
\max_{0\le n\le N_T} e_n^x
  \le
    C_T
    \left(
        e_0^x
        +
        \delta(\epsilon,N)
        +
        (\Delta t)^{p+1}
        +
        (\Delta t)^q
    \right)
,
    \]
which is the desired conclusion.

\section{Quantitative results for the Kuramoto-Sivashinsky example}
\label{app:KS_results}

Table~\ref{table:KSTrainingSizeMemory} reports the test valid prediction times (VPTs) for the KRR closure models with memory depths $m=0$ and $m=1$.
The $m=0$ results are summarized graphically in Figure~\ref{fig:KSTrainingSize}.
Across the training sizes considered here, adding one lag does not produce a systematic improvement for either the RBF or the diffusion maps kernel.
The same qualitative behavior was observed in the additional longer-memory experiments discussed in Section~\ref{sec:KS}.

Table~\ref{table:lstm_etdrk4} reports the VPTs for the LSTM closure models.
The first row corresponds to the result reported in \cite{harlim_2021}, employed on $2.5\times 10^5$ training samples and memory $m=19$.
In this case, a single model is trained and then evaluated on 50 test trajectories.

To obtain the results of the LSTM model, random Gaussian noise with standard deviation 0.1 is added to the training pairs, and the training pairs are normalized to have mean 0 and standard deviation 1.
Without this added noise, the VPTs are lower.
For each row, 10 models are trained on different trajectory segments and then evaluated on 50 test trajectories, resulting in a total of 500 tests for each case.
We note that as with the KRR models in Table~\ref{table:KSTrainingSizeMemory}, the change from $m=0$ to $m=1$ does not provide a significant systematic improvement over the range of training sizes considered.

\begin{table}[htbp]
    \centering
    \begin{tabular}{ccccc}
        \toprule
        $N_{\mathrm{train}}$
        &
        \multicolumn{2}{c}{$m=0$}
        &
        \multicolumn{2}{c}{$m=1$}
        \\
        \cmidrule(lr){2-3}
        \cmidrule(lr){4-5}
        &
        RBF
        &
        Diffusion
        &
        RBF
        &
        Diffusion
        \\
        \midrule
        2048
        &
        $1.0171\pm0.1177$
        &
        $0.9121\pm0.1183$
        &
        $0.9885\pm0.1562$
        &
        $0.9161\pm0.1567$
        \\
        4096
        &
        $1.1241\pm0.1884$
        &
        $1.1101\pm0.1780$
        &
        $1.2156\pm0.1276$
        &
        $1.1494\pm0.1948$
        \\
        8192
        &
        $1.6903\pm0.3016$
        &
        $1.8164\pm0.3150$
        &
        $1.6484\pm0.3067$
        &
        $1.8659\pm0.2142$
        \\
        16384
        &
        $2.3650\pm0.4693$
        &
        $\mathbf{2.6516\pm0.3819}$
        &
        $2.1400\pm0.4798$
        &
        $2.6490\pm0.4311$
        \\
        \bottomrule
    \end{tabular}
    \caption{
        Mean test VPT for the KRR closure models with memory depths $m=0$ and $m=1$.
        For each configuration, the reported statistics are the mean and standard deviation across ten independently trained-model test means, with $50$ test rollouts per trained model.
    }
    \label{table:KSTrainingSizeMemory}
\end{table}

\begin{table}[htbp]
    \centering
    \begin{tabular}{ccc}
        \toprule
        \textbf{$N_{\textnormal{train}}+2N_V$} & \textbf{Memory $m$} &
        \textbf{VPTs (ETDRK4)}
        \\
        \midrule
        $2.5 \times 10^5$ & 19 & $2.494688 \pm 0.95904$ \\
        \midrule
        2048+4000 & 0 & $1.41372390\pm0.6669579
        $ \\
        4096+4000 & 0 & $1.4588309 \pm 0.71275799
        $ \\
        8192+4000 & 0 & $1.5381616 \pm 0.725299612$ \\
        16384+4000 & 0 & $1.7413667 \pm 0.8070518383
        $ \\
        \midrule
        2048+4000 & 1 & $1.3993662\pm0.6346926$ \\
        4096+4000 & 1 & $1.4566895 \pm 0.6997771
        $ \\
        8192+4000 & 1 & $1.5757393 \pm 0.7621465$ \\
        16384+4000 & 1 & $1.7931989 \pm 0.82408019
        $ \\
        \bottomrule
    \end{tabular}
    \caption{
        Test VPTs of true 6 modes and predicted 6 modes for the LSTM model using ETDRK4.
    }
    \label{table:lstm_etdrk4}
\end{table}

\section{Kinetic solver, datasets, and experimental protocol}
\label{app:kinetic}

Appendix~\ref{app:kinetic} provides additional details for the kinetic experiments.
Section~\ref{app:kinetic_solver} describes the kinetic and macroscopic numerical solvers, Section~\ref{app:kinetic_data} specifies the training, validation, and test datasets together with the construction of the KRR regression data, and Section~\ref{app:kinetic_additional_figures} presents supplementary rollout figures.

\subsection{Kinetic solver}
\label{app:kinetic_solver}

In this Appendix, we describe the numerical solver used to generate reference trajectories for the experiments in Section~\ref{sec:kinetic}.

Let $f_{i,k}^n \approx f(x_i,v_k,t^n)$ denote the discrete distribution on a uniform periodic spatial mesh, $I_i=[x_{i-1/2},x_{i+1/2}]$, for $i=0,\ldots,N_x-1$, and a midpoint velocity grid $\{v_k\}_{k=0}^{N_v-1}$.
We use $N_x=256$, $N_v=128$, $x\in[-1,1]$, and $v\in[-8,8]$ so that $\Delta x=\frac{2}{256}=0.0078125$, and $\Delta v=\frac{16}{128}=0.125$.
The relaxation parameter is $\mathrm{Kn}=10^{-3}$.

With $\psi(v) = (1,v,\frac{1}{2}v^2)^\top$, the discrete conservative moments are
    \[
\bfU_i
  =
    \left\langle
        \psi f_i
    \right\rangle_h
  =
    \Delta v
    \sum_{k=0}^{N_v-1}
      \psi(v_k) f_{i,k}
,
    \]
where $\langle g\rangle_h = \Delta v\sum_{k=0}^{N_v-1}g_k, $ for a velocity-dependent quantity $g_k$.

\paragraph{Transport discretization.}

The transport term is discretized by a second-order finite-volume method \cite{PareschiRusso2005,habbershaw2025year4}.
Piecewise-linear interface values are obtained using a $\theta$-minmod slope limiter with limiter parameter $\theta=2$.
At each interface,
    \[
f_{i+1/2,k}^{-}
  =
    f_{i,k}+\frac{1}{2}s_{i,k}
,\qquad
f_{i+1/2,k}^{+}
  =
    f_{i+1,k}-\frac{1}{2}s_{i+1,k}
,
    \]
where $s_{i,k}$ denotes the limited reconstruction increment.
The velocity-upwind distribution flux is
    \[
\mathcal F_{i+1/2,k}
  =
    v_k^{+}f_{i+1/2,k}^{-}
    +
    v_k^{-}f_{i+1/2,k}^{+}
,\qquad
v_k^{+}
  =
    \max(v_k,0)
,\quad
v_k^{-}
  =
    \min(v_k,0)
.
    \]
The velocity moments of this flux give the numerical interface flux in \eqref{eq:kinetic_numerical_moment_flux}, which is used for model training.

\paragraph{IMEX time integration.}

Time integration uses the two-stage, second-order IMEX Runge--Kutta
method with coefficient matrices and vectors
    \[
\begin{array}{c|cc}
    0 & 0 & 0\\
    1 & 1 & 0\\
    \hline
    & \frac{1}{2} & \frac{1}{2}
\end{array}
\qquad\text{and}\qquad
\begin{array}{c|cc}
    \gamma & \gamma & 0\\
    1-\gamma & 1-2\gamma & \gamma\\
    \hline
    & \frac{1}{2} & \frac{1}{2}
\end{array}
,\quad
\gamma
  =
    1-\frac{1}{\sqrt2}
    \]
for the explicit transport and implicit relaxation Butcher tableaux, respectively.
Thus, transport is evaluated explicitly at each Runge--Kutta stage, while the BGK relaxation term is treated implicitly.

\paragraph{Conservative discrete Maxwellian.}

The continuum Maxwellian in \eqref{eq:bgk_maxwellian}, when simply sampled on a finite velocity grid, does not in general have quadrature moments exactly equal to the prescribed discrete moments.
To prevent the resulting collision-moment defect from entering the reference trajectories, the discrete kinetic solver uses the conservative discrete Maxwellian of Mieussens~\cite{Mieussens2000}.

For a target conservative state $\bfU = (\rho, \rho u, E)^\top$ the discrete equilibrium is written as
    \[
M_{h,k}(\boldsymbol{\alpha})
  =
    \exp
    \left(
        \alpha_0
        +
        \alpha_1 v_k
        +
        \alpha_2\frac{v_k^2}{2}
    \right)
,
    \]
where the parameters $\boldsymbol{\alpha}=(\alpha_0,\alpha_1,\alpha_2)^\top$ are chosen so that
    \begin{equation}
\Delta v
\sum_{k=0}^{N_v-1}
  \psi(v_k)
M_{h,k}(\boldsymbol{\alpha})
  =
    \bfU
.
    \label{eq:discrete_maxwellian_constraints}
    \end{equation}
Specifically, after eliminating $\alpha_0$ through the mass constraint, the implementation solves a two-variable convex optimization problem for $(\alpha_1,\alpha_2)$ using a damped Newton iteration.
The solution is accepted according to the residual in the full three physical moment equations \eqref{eq:discrete_maxwellian_constraints}.

Writing $\boldsymbol{\beta}=(\alpha_1,\alpha_2)^\top$, $\boldsymbol{z}_k=(v_k,v_k^2/2)^\top$, and
    \[
Z(\boldsymbol{\beta})
  =
    \sum_k \omega_k
    \exp(\boldsymbol{\beta}^\top\boldsymbol{z}_k)
,
    \]
the mass constraint gives $\alpha_0=\log\rho-\log Z(\boldsymbol{\beta})$.
The remaining parameters minimize the convex dual problem
    \[
\Phi(\boldsymbol{\beta})
  =
    \log Z(\boldsymbol{\beta})
    -
    \boldsymbol{\beta}^\top
    \begin{pmatrix}
        u\\ E/\rho
    \end{pmatrix}
.
    \]
We solve this two-dimensional problem using a damped Newton iteration and accept the solution according to the residual in the full three moment constraints \eqref{eq:discrete_maxwellian_constraints}.

The same conservative discrete-equilibrium construction is used both to initialize the kinetic distribution from prescribed macroscopic fields and to evaluate the equilibrium in the BGK collision term.
Consequently, the collision contribution has vanishing discrete mass, momentum, and energy moments up to the tolerance of the nonlinear solve.

\paragraph{Time step and saved trajectory.}

The fixed time step is $\Delta t = 0.9\,\frac{\Delta x}{2\max_k |v_k|} = 4.39453125\times10^{-4}.$
Each standard trajectory contains $455$ time intervals and $456$ saved time levels.
The final saved time is therefore $t_{\mathrm F} = 455\Delta t = 0.199951171875$.

\paragraph{Macroscopic solver.}

The independently evolved macroscopic model uses the same spatial mesh and the same $\theta$-minmod piecewise-linear reconstruction procedure as the kinetic solver, applied componentwise to the conservative macroscopic variables rather than to the discrete distribution function.
The reconstructed left and right macroscopic states are then used to construct the zero-closure Rusanov (local Lax--Friedrichs) numerical flux \cite{Rusanov1962,Toro2009}, as defined in \eqref{eq:macroscopic_numerical_flux}, with local wave speed $a_{i+1/2} = \max\left\{ |u^-|+\sqrt{3\theta^-}, |u^+|+\sqrt{3\theta^+} \right\}$.
The kinetic and macroscopic solvers therefore construct their transport fluxes independently.
The learned target $\bfR_{\mathrm{eff}}$ in \eqref{eq:R_eff_definition} represents the correction between the corresponding stage-weighted numerical transport fluxes.

\subsection{Kinetic datasets and construction of regression data}
\label{app:kinetic_data}

All reference trajectories are initialized from the macroscopic profiles in \eqref{eq:macroscopic_ic_fields}, with the corresponding kinetic distribution initialized using the conservative discrete equilibrium described in Section~\ref{app:kinetic_solver}.
We organize the initial conditions into two data types: Type I data contains unimodal profiles ($P=1$), while Type II data contains bimodal profiles ($P=2$).

\paragraph{Training, validation, and testing dataset families.}

The Type I training data set consists of the $3\times3\times8$ product of combinations of
    \[
A
  \in
    \{1.5,2.0,2.5\}
,\quad
\sigma
  \in
    \{0.04,0.05,0.06\}
,\quad
\eta
  \in
    \left\{
        -\frac{1}{2},
        -\frac{3}{8},
        -\frac{1}{4},
        -\frac{1}{8},
        0,
        \frac{1}{8},
        \frac{1}{4},
        \frac{3}{8}
    \right\}
,\quad
\mu
  =
    \eta\Delta x
.
    \]
This gives $72$ unimodal training trajectories.

The Type II training data set contains $72$ bimodal initial conditions
parameterized by
    \[
A_1,A_2
  \in
    [1.4,2.5]
,\qquad
\sigma_1,\sigma_2
  \in
    [0.04,0.06]
,\qquad
\mu_1
  =
    \phi-\frac{d}{2}
,\qquad
\mu_2
  =
    \phi+\frac{d}{2}
,
    \]
where $d\in[0.35,0.90]$, and $\phi\in[-0.375,0.375]$.
The $72$ parameter sets are selected using a fixed Latin-hypercube design with an additional space-filling optimization over this six-dimensional parameter domain \cite{McKayBeckmanConover1979,MorrisMitchell1995}.

Hyperparameter selection uses a common set of six Type I and six Type II validation trajectories for all three closure representations.
These validation cases are distinct from the training cases and fixed before their trajectory outcomes were generated.

Final model testing uses $24$ Type I and $24$ Type II initial conditions that are disjoint from the training and validation sets.
These testing trajectories are not used during preprocessing, hyperparameter selection, or model development.

The exact initial-condition parameters for all training, validation, and test trajectories, together with the scripts used to construct them, are provided in the public repository.

\paragraph{Regression targets.}

For each time interval $[t^n,t^{n+1}]$, the kinetic and zero-closure macroscopic solvers are advanced independently through their Runge--Kutta stages.
Their stage-weighted numerical transport fluxes define the supervised target
    \[
\bfR_{\mathrm{eff},i+1/2}^n
  =
    \mathcal{G}_{\mathrm{kin},i+1/2}^n
    -
    \mathcal{G}_{\mathrm{mac},i+1/2}^n
,
    \]
as in \eqref{eq:R_eff_definition}.
Thus the regression target is the three-component numerical-flux correction rather than the physical heat flux directly.

\paragraph{Construction of KRR training sets.}

For each global and local model, we compare the use of Type I training data with using mixed Type I+II training data, while holding the total number of KRR training pairs fixed.
For the global model, each training sample corresponds to a complete spatial-field snapshot in $\bbR^{3N_x}$.
The model trained on Type I data uses $114$ approximately uniformly sampled $3N_x$--dimensional snapshots from each of the $72$ Type I training trajectories, while the mixed model uses $57$ snapshots from each of the $72$ Type I and $72$ Type II training trajectories.

For the local models, training pairs are instead selected from a sparse subset of the available space--time interface data spanning the temporal evolution.
The zeroth-order model trained on Type I data uses $56$ local interface samples from each of the $72$ Type I training trajectories, while the mixed model uses $28$ samples from each of the $72$ Type I and $72$ Type II training trajectories.
The first-order residual model follows the same general allocation strategy with twice as many training samples.

The resulting allocations are summarized in Table~\ref{table:kinetic_center_budgets}.
The exact sample indices and data point identities are provided in the public repository.

\begin{table}[H]
    \centering
    \begin{tabular}{lcc}
    \toprule
    Representation & Type I & Type I+II (Mixed)
    \\
    \midrule
    Global
    & $72 \times 114 = 8208$
    & $72 \times 57 +72 \times 57=8208$
    \\
    Zeroth-order local
    & $72 \times 56=4032$
    & $72 \times 28 + 72 \times 28 = 4032$
    \\
    First-order local
    & $72 \times 112 = 8064$
    & $72 \times 56 + 72 \times 56 = 8064$
    \\
    \bottomrule
    \end{tabular}
    \caption{
        KRR training-pair allocations used in the kinetic-model comparisons.
        In each table entry, the number $72$ represents the number of Type I and Type II trajectories used for training.
        For the global model, the second factor denotes the number of temporal full-field snapshots selected per trajectory; for the local models, it denotes the number of local space--time interface samples selected per trajectory.
    }
    \label{table:kinetic_center_budgets}
\end{table}

For the first-order model, a single Type I zeroth-order model with $N_0=4032$ centers is trained first and then held fixed.
The residual correction is trained using the derivative representation defined in Section~\ref{sec:kinetic_representations}, and the same baseline is used in both first-order comparisons.

See the public repository for the exact temporal and spatial sample indices, local-activity ranking and ordering rules, normalization statistics, duplicate-filtering convention, and final KRR center identities used in each experiment.

\subsection{Additional kinetic rollout figures}
\label{app:kinetic_additional_figures}

\begin{figure}[htbp]
    \centering
    \includegraphics[
        width=0.99\linewidth
    ]{
        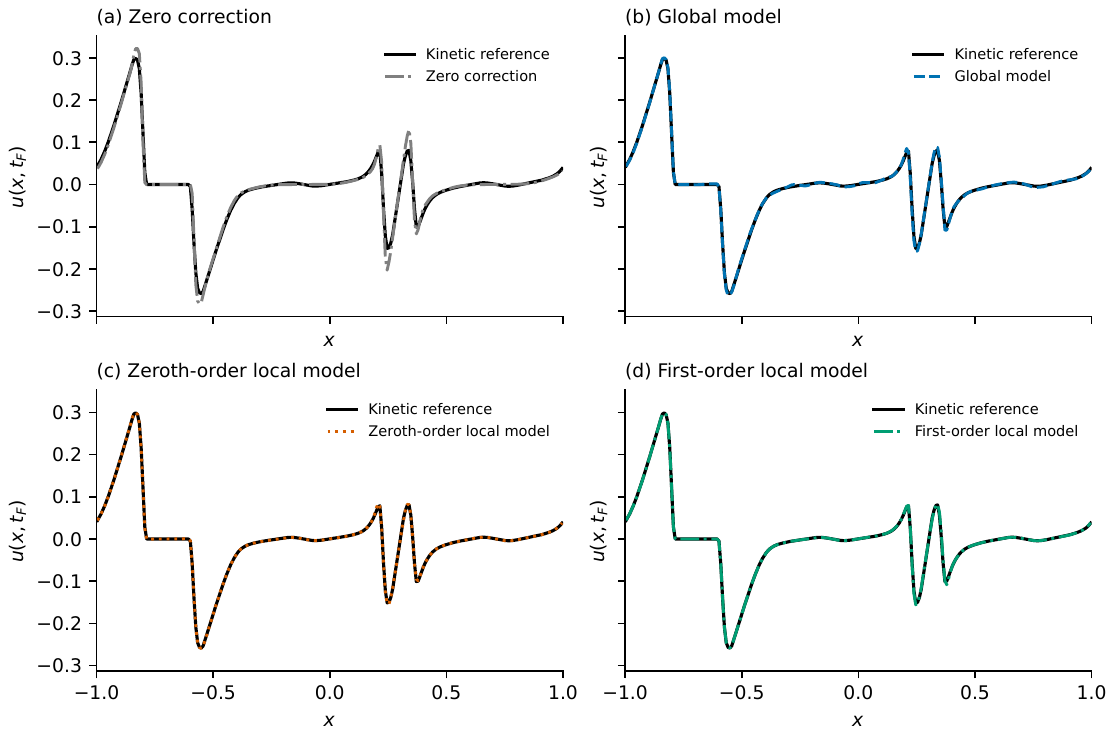
    }
    \caption{
        Final-time velocity profiles for the same bimodal rollout used in the main-text kinetic comparison in Figure~\ref{fig:kinetic_pairwise_rho}.
        Each panel compares the kinetic reference with one macroscopic model: the zero-correction baseline, the global closure, the zeroth-order local closure, and the first-order local residual closure.
    }
    \label{fig:kinetic_additional_u}
\end{figure}

\begin{figure}[htbp]
    \centering
    \includegraphics[
        width=0.99\linewidth
    ]{
        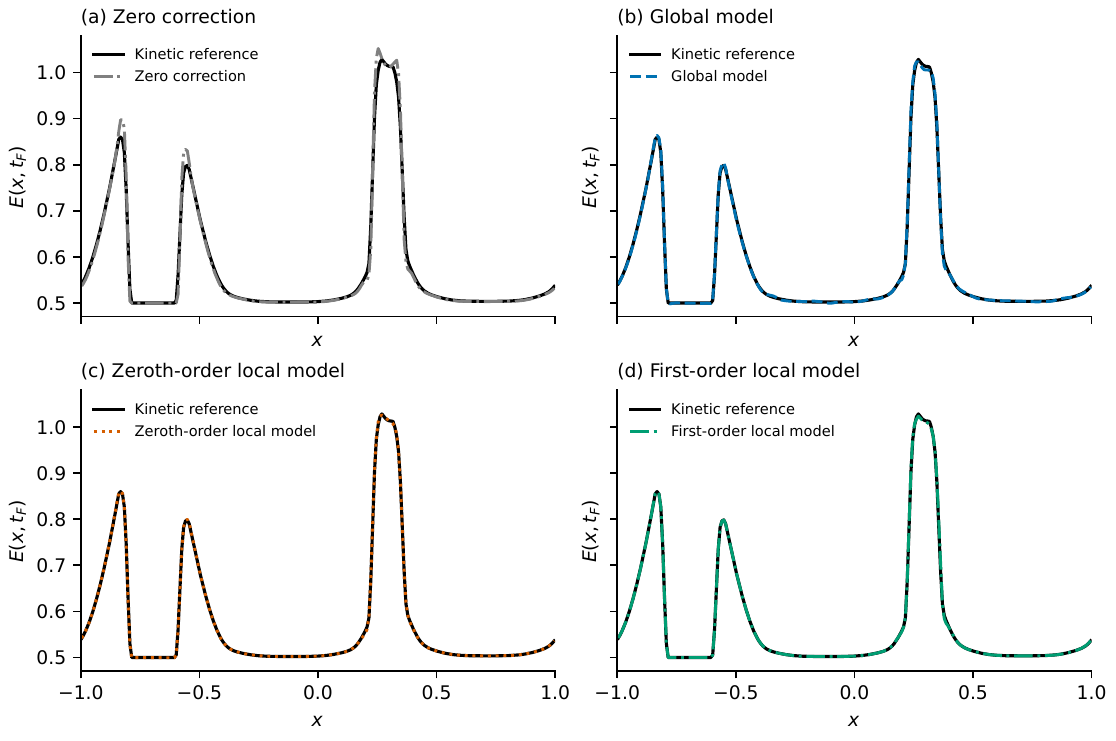
    }
    \caption{
        Final-time energy profiles for the same bimodal rollout used in the main-text kinetic comparison in Figure~\ref{fig:kinetic_pairwise_rho}.
        Each panel compares the kinetic reference with one macroscopic model, using the same model ordering as Figure~\ref{fig:kinetic_additional_u}.
    }
    \label{fig:kinetic_additional_E}
\end{figure}

Figures~\ref{fig:kinetic_additional_u}--\ref{fig:kinetic_additional_E} provide supplementary final-time comparisons of velocity, temperature, and energy, respectively, in the same setting as Figure~\ref{fig:kinetic_pairwise_rho}.

\begin{figure}[htbp]
    \centering
    \includegraphics[
        width=0.99\linewidth
    ]{
        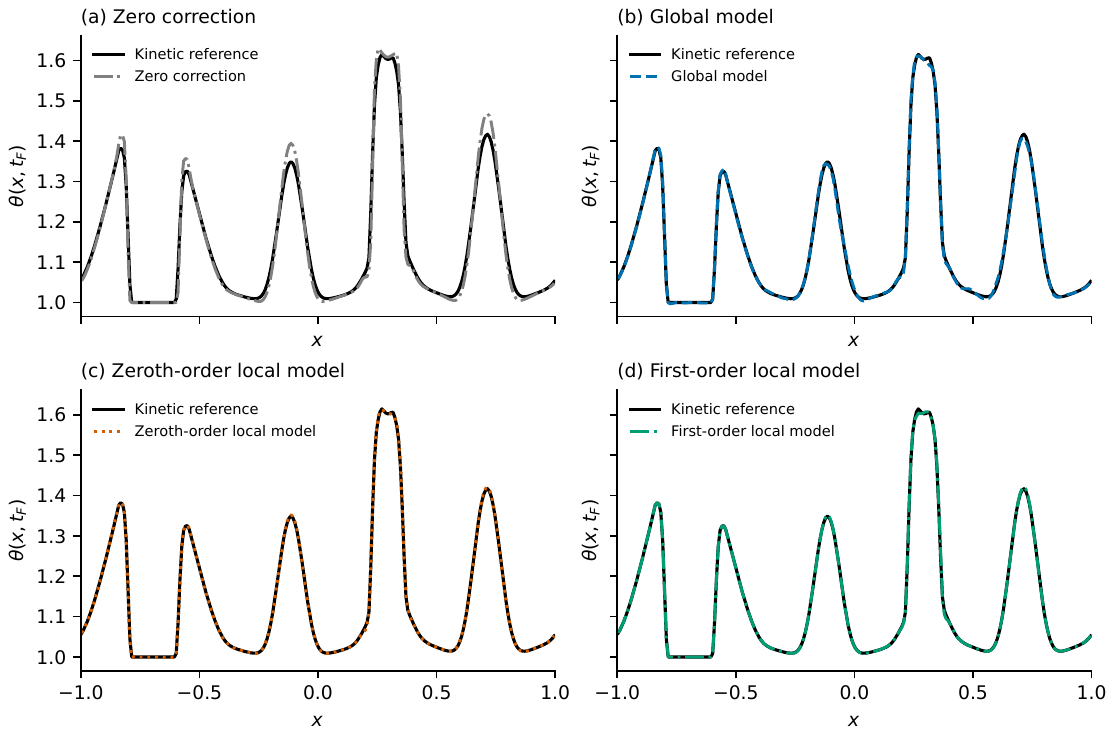
    }
    \caption{
        Final-time temperature profiles for the same bimodal rollout used in the main-text kinetic comparison in Figure~\ref{fig:kinetic_pairwise_rho}.
        Each panel compares the kinetic reference with one macroscopic model, using the same model ordering as Figure~\ref{fig:kinetic_additional_u}.
    }
\end{figure}

\end{document}